# MAXIMUM AND ENTROPIC REPULSION FOR A GAUSSIAN MEMBRANE MODEL IN THE CRITICAL DIMENSION[1]

By Noemi Kurt

*Universität Zürich*

We consider the real-valued centered Gaussian field on the four-dimensional integer lattice, whose covariance matrix is given by the Green's function of the discrete Bilaplacian. This is interpreted as a model for a semiflexible membrane. $d = 4$ is the critical dimension for this model. We discuss the effect of a hard wall on the membrane, via a multiscale analysis of the maximum of the field. We use analytic and probabilistic tools to describe the correlation structure of the field.

**1. Introduction and main results.** Let $V := [-1, 1]^d$, and $V_N := NV \cap \mathbb{Z}^d$. In this paper we consider the real-valued Gaussian field $\varphi = \{\varphi_x\}_{x \in V_N}$, whose covariance matrix is given by the Green's function of the discrete Bilaplacian. Such a field can be interpreted as a model for a $d$-dimensional interface in $d + 1$-dimensional space. It is described by the formal Hamiltonian $H_N(\varphi) = \frac{1}{2} \sum_x (\Delta \varphi_x)^2$. For this model, $d = 4$ is critical in the sense that, in dimensions higher than 4, the infinite volume Gibbs measure exists (see [10, 13]), but not in $d = 4$ and below. A phenomenon of interest for random interface models is the so-called entropic repulsion, which refers to the fact that the presence of a hard wall forces the interface to move away from the wall. This is modeled by requiring the field $\{\varphi_x\}$ to be positive inside a certain region. To mathematically understand entropic repulsion, one needs to study the asymptotics of the probability $P(\varphi_x \geq 0, x \in V)$ for some region $V \subset \mathbb{Z}^d$. In the case considered in this paper, this is achieved by first investigating the asymptotic behavior of the maximum of the field, via a sophisticated multiscale-analysis developed in [1] for the lattice free field in the critical dimension. The main difficulty is due to the fact that, unlike the lattice free

Received January 2008; revised June 2008.

[1]Supported in part by the Swiss National Science Foundation Contract 200020-116348.

field, our model does not have a random walk representation, which is crucial in most approaches to the lattice free field (see, e.g., [1, 2]). To obtain the analogous results, we use methods from PDE to get good estimates of some discrete biharmonic Green's functions.

For $k \in \mathbb{N}$, let $\partial_k V_N := \{x \in V_N^c : \text{dist}(x, V_N) \leq k\}$ be the boundary of thickness $k$ of $V_N$. We write $\partial V_N := \partial_1 V_N$ for the simple boundary. The discrete Laplacian $\Delta$ is defined on functions $f : \mathbb{Z}^d \to \mathbb{R}$ by

$$\Delta f(x) := \frac{1}{2d} \sum_{i=1}^{d} (f(x + e_i) + f(x - e_i) - 2f(x)),$$

where $e_i$ denotes the unit vector in the $i$th coordinate direction. With some abuse of notation, we write $\Delta f_x := (\Delta f)(x)$. By $\Delta_N$, we denote the restriction of this operator to functions which are equal to 0 outside $V_N$. We write $\Delta^2$ for the iteration, $\Delta^2 f(x) := \Delta(\Delta f)(x)$, and $\Delta_N^2$ for the restriction of $\Delta^2$ to functions which are equal to 0 outside $V_N$. It is important to notice that $\Delta_N^2 \neq (\Delta_N)^2$. We can view $\Delta_N^2$ as the matrix given by

$$\Delta_N^2(x, y) = \begin{cases} 1 + \dfrac{1}{2d}, & \text{if } x = y, x \in V_N, \\[2mm] -\dfrac{1}{d}, & \text{if } |x - y| = 1, x, y \in V_N, \\[2mm] \dfrac{1}{4d^2}, & \text{if } |x - y| = 2, x, y \in V_N, \\[2mm] \dfrac{1}{2d^2}, & \text{if } |x - y| = \sqrt{2}, x, y \in V_N, \\[2mm] 0, & \text{otherwise.} \end{cases}$$

The matrix $(\Delta_N^2(x, y))_{x, y \in V_N}$ is positive definite (see Remark A.7). Let $G_N(x, y)$ be it's matrix inverse. This means that we can interpret $G_N$ as a Green's function given by the following discrete biharmonic boundary value problem on $V_N$: For $x \in V_N$,

$$\begin{aligned} \Delta^2 G_N(x, y) &= \delta(x, y), & y \in V_N, \\ G_N(x, y) &= 0, & y \in \partial_2 V_N. \end{aligned} \tag{1}$$

To see the connection to boundary value problems of PDE, note that this is a discrete version of the (continuous) biharmonic boundary value problem with Dirichlet boundary conditions:

$$\begin{aligned} \Delta^2 u(x) &= f(x), & x \in V, \\ u(x) &= 0, & x \in \partial V, \\ \frac{d}{dn} u(x) &= 0, & x \in \partial V. \end{aligned}$$

Here, $\frac{d}{dn}$ denotes the derivative in the direction of the outer normal vector. However, we will not directly use this correspondence between discrete and continuous, apart from gaining inspiration from standard PDE methods.

The model we study in this paper is the centered Gaussian field $\{\varphi_x\}_{x \in V_N}$ on $V_N$ with covariances $\mathrm{cov}_N(\varphi_x, \varphi_y) = G_N(x, y)$. Denote the law of this field by $P_N$. Algebraic manipulations show that $P_N$ is the Gibbs measure on $\mathbb{R}^{V_N}$ with 0 boundary conditions outside $V_N$ and Hamiltonian

$$H_N(\varphi) = \frac{1}{2} \sum_{x \in \mathbb{Z}^d} (\Delta \varphi_x)^2.$$

Note that the choice of boundary conditions in the definition of $\Delta_N^2$ and $G_N$ is absolutely crucial in order to obtain a Gibbs measure (see [6], Chapter 13), meaning that, for $A \subset V_N$, the distribution conditional on $\mathcal{F}_{A^c} = \sigma(\varphi_x, x \in A^c)$, the sigma field generated by $\varphi_x, x \in A^c$, satisfies

$$P_N(\cdot | \mathcal{F}_{A^c})(\psi) = P_{A,\psi}(\cdot), \qquad P_N(d\psi)\text{-a.s.},$$

where

$$P_{A,\psi}(d\varphi) := \frac{1}{Z_A} \exp\left( -\frac{1}{2} \sum_{x \in \mathbb{Z}^d} (\Delta \varphi_x)^2 \right) \prod_{x \in A} d\varphi_x \prod_{x \in V_N \setminus A} \delta_{\psi_x}(d\varphi_x).$$

($Z_A$ is the normalizing constant.) This implies that $P_N(\cdot | \mathcal{F}_{A^c})$ is the Gaussian distribution with mean

$$(2) \qquad m_x = -\sum_{y \in A} (\Delta_A^2)^{-1}(x, y) \sum_{z \in A^c} \Delta^2(y, z) \psi_z$$

and covariance matrix $(\Delta_A^2)^{-1}$. Here, $\Delta_A^2$ is the restriction of $\Delta^2$ to functions, which are 0 outside $A$. We would not obtain this Gibbsianness if we chose $(\Delta_N)^2$ [resp. $(\Delta_A)^2$] in the place of $\Delta_N^2$ (resp. $\Delta_A^2$). Since the range of interaction of $\Delta^2$ is 2, we see that $P_N(\cdot | \mathcal{F}_{A^c}) = P_N(\cdot | \mathcal{F}_{\partial_2 A})$.

This model is called the membrane or Laplacian model. One should compare it to the well-known lattice free field or gradient model, whose Hamiltonian is given by $H_N^\nabla(\varphi) := \frac{1}{2d} \sum |\nabla \varphi_x|^2$. Note that $H_N^\nabla(\varphi)$ is small if $\varphi$ is approximately constant, which implies that this model favors interfaces that are essentially flat. On the other hand, the membrane model prefers configurations with constant curvature. In the physics literature, for example [9, 14], linear combinations of the two models are considered as models for semiflexible membranes (or semiflexible polymers if $d = 1$). Contrary to the gradient model, there are only a few mathematically rigorous results for the membrane model, in $d \geq 5$, where the infinite volume limit exists [10, 13], and in $d = 1$ [3, 4]. One reason why the Laplacian model is more difficult to study is the absence of a random walk representation, which is exploited for the gradient model, and allows to get precise expressions for many quantities, in particular, the variance. In this paper we treat the membrane model in the critical dimension, $d = 4$, which means that we need to consider the finite volume $V_N$, where boundary effects come into play. Although we do not investigate the behavior of the field close to the boundary but only in the bulk, there are considerable analytical difficulties to overcome, which stem from the boundary conditions of the Green's function. We are able, using analytical and probabilistic methods, to control the variances in a way that is sufficient to apply the methods of [1]. Let $\gamma := \frac{8}{\pi^2}$, and define for $\delta \in (0, 1/2)$

$$V_N^\delta := \{x \in V_N : \operatorname{dist}(x, V_N^c) \geq \delta N\}.$$

Our first result consists of bounds on the variances:

PROPOSITION 1.1. *Let $d = 4$, and let $0 < \delta < 1/2$.*

(a) *There exists $C > 0$ such that $\sup_{x \in V_N} \operatorname{var}_N(\varphi_x) \leq \gamma \log N + C$.*

(b) *There exists $C(\delta) > 0$ such that $\sup_{x \in V_N^\delta} |\operatorname{var}_N(\varphi_x) - \gamma \log N| \leq C(\delta)$.*

Proposition 1.1 (together with the concentration result Lemma 2.11 in the next section) is the key to the results in this paper. It shows why the four-dimensional membrane model behaves in many ways like the two-dimensional lattice free field. We have the same behavior of the maximum:

THEOREM 1.2. *Let $d = 4$.*

(a)

$$\lim_{N \to \infty} P_N\left(\sup_{x \in V_N} \varphi_x \geq 2\sqrt{2\gamma} \log N\right) = 0$$

(b) *Let $0 < \delta < 1/2$, and $0 < \eta < 1$. There exists a constant $c = c(\eta, \delta) > 0$, such that*

$$P_N\left(\sup_{x \in V_N^\delta} \varphi_x \leq (2\sqrt{2\gamma} - \eta) \log N\right) \leq \exp(-c(\log N)^2).$$

These bounds on the maximum allow us to give the precise asymptotics of the probability that the field is positive on a certain region inside $V_N$. Let $D \subset V$ be connected with smooth boundary, which has positive distance to $\partial V$. Let $D_N := ND \cap \mathbb{Z}^4$ and define

$$\Omega_N^+ := \{\{\varphi_x\}_{x \in V_N} : \varphi_x \geq 0 \ \forall x \in D_N\}.$$

We think of $D_N$ as a hard wall that forces the field to be positive. The probability of this event is given by our next result. Let $H^2(V)$ denote the usual Sobolev space of twice differentiable functions on $V$, and $H_0^2$ the subspace of functions in $H^2(V)$ which are 0 at the boundary of $V$.

THEOREM 1.3.  *Let $d = 4$.*

$$\lim_{N \to \infty} \frac{1}{(\log N)^2} \log P_N(\Omega_N^+) = -8\gamma \mathcal{C}_V^2(D),$$

*where $\mathcal{C}_V^2(D) = \inf\{\frac{1}{2} \int_V |\Delta h|^2 \, dx : h \in H_0^2(V), h \geq 1 \text{ a.e. on } D\}$.*

One would like to understand the behavior of the field conditioned on the event $\Omega_N^+$. We can prove the following. For $0 \leq \varepsilon < 1$ and $x \in D_N$, let $V_{\varepsilon N}(x)$ denote the box of side-length $\varepsilon N$ with center $x$, and $\overline{\varphi}_{\varepsilon N}(x) := \frac{1}{|V_{\varepsilon N}(x)|} \sum_{y \in V_{\varepsilon N}(x)} \varphi_y$.

PROPOSITION 1.4.  *For any $\eta > 0$,*

$$\lim_{N \to \infty} \sup_{\substack{x \in D_N, \\ V_{\varepsilon N}(x) \subset D_N}} P_N(\overline{\varphi}_{\varepsilon N}(x) \leq (2\sqrt{2\gamma} - \eta) \log N | \Omega_N^+) = 0.$$

This implies that the local sample mean of the field is pushed by the hard wall to a height of at least $2\sqrt{2\gamma} \log N$. In the physics literature this phenomenon is referred to as entropic repulsion [12], since it is due to the fluctuations of the field that it moves away from the wall. It is expected that the upper bound on the height of the conditioned field is the same, that is, that $P_N(\overline{\varphi}_{\varepsilon N}(x) \geq (2\sqrt{2\gamma} + \eta) \log N | \Omega_N^+) = 0$. Also, for the gradient model, the result holds for the height variables $\varphi_x$ in the place of $\overline{\varphi}_{\varepsilon N}(x)$ [1]. The proof for the gradient model uses the FKG-inequalities. For the membrane model, the criterion for the FKG-property, Corollary 1.8 of [8] is satisfied only in the infinite volume case and without the positivity constraint. We therefore need to average over the heights in order to obtain the result.

The paper is organized as follows. In the next section we investigate the variance structure of the four-dimensional membrane model and prove Proposition 1.1 and some related results. Here we exploit the fact that we can compare $G_N$ to the Green's function corresponding to $(\Delta_N)^2$, for which we have a random walk interpretation. The comparison of the two Green's functions is based on analytical tools on the regularity of the solutions of boundary value problems. Some of the more technical proofs are deferred to the Appendix. In Section 3 we give the proof of Theorem 1.2, using the same multiscale analysis as for the gradient model. We refer to [1] for detailed comments on the ideas behind this method. The proof of Theorem 1.3 is given in Section 4, and that of Proposition 1.4 in Section 5.

Throughout the paper, $c, C, c'$ etc. will denote generic positive constants whose value may change from line to line. By $B_r$ we denote the ball of radius $r$ and center 0.

**2. Variance structure and the discrete Green's function.** The aim of this section is to control $G_N(x,x)$. To this purpose, we compare it to a biharmonic Green's function with different boundary conditions. Let

$$E_1 := \{v : V_N \cup \partial_2 V_N \to \mathbb{R} : v(x) = 0 \ \forall x \in \partial_2 V_N\}.$$

Recall from the Introduction that the covariance matrix of the model is given by the unique function $G_N(x,\cdot)$ in $E_1$ which satisfies $\Delta^2 G_N(x,y) = \delta(x,y)$.

Let us introduce the usual harmonic Green's function. Let $A$ be an arbitrary subset of $\mathbb{Z}^d$, fix $x \in A$, and let $\Gamma_A(x,\cdot)$ be the unique lattice function which satisfies

$$\Delta \Gamma_A(x,y) = -\delta(x,y), \qquad y \in A,$$
$$\Gamma_A(x,y) = 0, \qquad y \in \partial A.$$

(Existence and uniqueness follows from standard discrete harmonic analysis; see, e.g., Chapter I of [11].) Let $\Gamma_N(x,y) := \Gamma_{V_N}(x,y)$. Define now for $x, y \in V_N$,

$$(3) \qquad \overline{G}_N(x,y) := \sum_{z \in V_N} \Gamma_N(x,z) \Gamma_N(z,y),$$

and extend $\overline{G}_N(x,\cdot)$ to a function on $V_N \cup \partial_2 V_N$ by requiring

$$(4) \qquad \begin{aligned} \overline{G}_N(x,y) &= 0, \qquad y \in V_{N+1} \setminus V_N \quad \text{and} \\ \Delta \overline{G}_N(x,y) &= 0, \qquad y \in \partial V_N. \end{aligned}$$

It is straightforward to check that, with these conditions, $\Delta^2 \overline{G}_N(x,y) = \delta(x,y)$ for all $x, y \in V_N$. In fact, $\overline{G}_N(x,\cdot)$ is the (again unique) function which satisfies

$$\Delta^2 \overline{G}_N(x,y) = \delta(x,y), \qquad y \in V_N,$$
$$\overline{G}_N(x,y) = 0, \qquad y \in V_{N+1} \setminus V_N,$$
$$\Delta \overline{G}_N(x,y) = 0, \qquad y \in \partial V_N.$$

The main idea of this section is to compare $G_N(x,y)$ and $\overline{G}_N(x,y)$. In fact, we will later on show that if $x \in V_N^\delta$,

$$\sup_{y \in V_N^\delta} |G_N(x,y) - \overline{G}_N(x,y)| \le c$$

for some $c = c(\delta) < \infty$. This will be done by studying the boundary value problem satisfied by $G_N(x,y) - \overline{G}_N(x,y)$ and showing that the solution of this boundary value problem is sufficiently regular (in a sense to be specified). Since $\overline{G}_N$ is given in terms of $\Gamma_N$, well-known results from harmonic analysis and random walks give us a very good control on the behavior of $\overline{G}_N(x, y)$. Combining all this will then prove Proposition 1.1.

Before embarking on the comparison of $G_N$ and $\overline{G}_N$, we derive the necessary estimates on $\overline{G}_N$. We collect the following well-known results on $\Gamma_N$, which we will use to describe $\overline{G}_N$. For proofs we refer to [11], Chapter I. Let $A$ be an arbitrary subset of $\mathbb{Z}^d$, and write $\Gamma_A$ for the Green's function of the Dirichlet problem on $A$. The following hold:

- $\Gamma_A(x, y)$ is the expected number of visits in $y \in A$ of a simple random walk starting at $x$ which is killed as it exits $A$, that is,

$$\Gamma_A(x, y) = \mathbb{E}^x\left(\sum_{k=0}^{\tau_A} 1_{\{X_k = y\}}\right) = \sum_{k=0}^{\infty} \mathbb{P}^x(X_k = y, k < \tau_A),$$

  where $\tau_A = \inf\{k \geq 0 : X_k \in A^c\}$.
- If $d \geq 3$, $\lim_{N \to \infty} \Gamma_{V_N}(x, y) =: \Gamma(x, y)$ exists for all $x, y \in \mathbb{Z}^d$, and as $|x - y| \to \infty$,

$$\Gamma(x, y) = a_d \frac{1}{|x - y|^{d-2}} + O(|x - y|^{1-d}),$$

  with $a_d = \frac{2}{(d-2)\omega_d}$, where $\omega_d$ is the volume of the unit ball in $\mathbb{R}^d$.
- ([11], Proposition 1.5.9) If $d \geq 3$, for all $x \neq 0$

$$\Gamma_{B_N}(0, x) = a_d\left(\frac{1}{|x|^{d-2}} - \frac{1}{N^{d-2}}\right) + O(|x|^{1-d}).$$

- If $d \geq 3$, then

$$\Gamma_A(x, y) = \Gamma(x, y) - \sum_{z \in \partial A} \mathbb{P}^x(X_{\tau_A} = z)\Gamma(z, y).$$

- $\Gamma_A(x, y) = \Gamma_A(y, x)$.
- $\Gamma_A(x, y) \leq \Gamma_B(x, y)$ if $A \subset B$.

The fact that $\overline{G}_N$ is just the convolution of $\Gamma_N$ with itself leads to the following representation in terms of simple random walk: Letting $x, y \in V_N$, let $\{X_k\}, \{Y_m\}$ be two independent simple random walks on the lattice $\mathbb{Z}^d$, whose joint law with start in $x$ and $y$ respectively we denote by $\mathbb{P}^{x,y}$. Let $\tau_N$ denote the first exit time of $V_N$. Now we see from the random walk representation of $\Gamma_N$ that

$$\overline{G}_N(x, y) = \sum_{z \in V_N} \Gamma_N(x, z)\Gamma_N(z, y) = \mathbb{E}^{x,y}\left[\sum_{k=0}^{\tau_N} \sum_{m=0}^{\tau_N} 1_{\{X_k = Y_m\}}\right]$$

and

$$\overline{G}_N(x, y) = \sum_{z \in V_N} \Gamma_N(x, z)\Gamma_N(z, y)$$

$$= \sum_{k,m=0}^{\infty} \sum_{z \in V_N} \mathbb{P}^x(X_k = z, k < \tau_N) \mathbb{P}^z(Y_m = y, m < \tau_N)$$

$$= \sum_{k,m=0}^{\infty} \mathbb{P}^x(X_{k+m} = y, k + m < \tau_N)$$

$$= \sum_{k=0}^{\infty} (k+1) \mathbb{P}^x(X_k = y, k < \tau_N).$$

Hence, we have proven the following:

LEMMA 2.1.    *If $x, y \in V_N$, the following hold:*

$$\overline{G}_N(x,y) = \mathbb{E}^{x,y} \left[ \sum_{k=0}^{\tau_N} \sum_{m=0}^{\tau_N} 1_{\{X_k = Y_m\}} \right] = \sum_{k=0}^{\infty} (k+1) \mathbb{P}^x(X_k = y, k < \tau_N).$$

Estimates on $\overline{G}_N(x,x)$ are easily obtained:

LEMMA 2.2.    *Let $d = 4$. If $\delta \in (0, 1/2)$, there exist constants $c_1 = c_1(\delta) > 0$, $c_2 > 0$, such that, for $x \in V_N^\delta$,*

$$\frac{8}{\pi^2} \log N + c_1 \leq \overline{G}_N(x,x) \leq \frac{8}{\pi^2} \log N + c_2.$$

PROOF.    Let $B_r(x)$ denote the ball of radius $r$ about $x \in V_N$. Since $\Gamma_N(x,x) \leq \Gamma(x,x)$, we obtain

$$\overline{G}_N(x,x) \leq \sum_{z \in B_{2N}} \Gamma(x,z)\Gamma(z,x) \leq a_4^2 \sum_{\substack{z \in B_{2N}(x) \\ z \neq x}} \frac{1}{|x-z|^4} + O(1)$$

$$\leq 4a_4^2 \omega_4 \int_1^{2N} \frac{1}{r} \, dr + O(1) = \frac{8}{\pi^2} \log(2N) + c.$$

The lower bound follows by taking $B_{\delta N}(x)$ in the place of $B_{2N}(x)$:

$$\overline{G}_N(x,x) \geq \sum_{z \in B_{\delta N}} \Gamma_{B_{\delta N}}(x,z)\Gamma_{B_{\delta N}}(z,x) \geq 4a_4^2 \omega_4 \int_1^{\delta N} \frac{1}{r} \, dr + O(1)$$

$$= \frac{8}{\pi^2} \log(\delta N) + c.$$

$\square$

We need to introduce discrete Sobolev norms. Let $\partial_- V_N := \{x \in V_N : \operatorname{dist}(x, V_N^c) \leq 1\}$. We denote the first difference in the $i$th direction of a function $v : \mathbb{Z}^d \to \mathbb{R}$ by $\nabla_i v(x) := v(x + e_i) - v(x)$, and more general, for a multiindex $\alpha = (\alpha_1, \ldots, \alpha_d) \in \mathbb{N}^d$, write $\nabla^\alpha v(x) := \nabla_1^{\alpha_1} \cdots \nabla_d^{\alpha_d} v(x)$.

For $v : V_N \cup \partial_k V_N \to \mathbb{R}$ define

$$\|v\|_{H^k(V_N)}^2 := \sum_{j=0}^k \sum_{\substack{\alpha \in \mathbb{N}^d: \\ |\alpha| = j}} \sum_{x \in V_N} (N^j \nabla^\alpha v(x))^2.$$

For $v, w \in E_1$ define

$$\mathcal{D}(v, w) := \sum_{x \in V_N} \Delta v(x) \Delta w(x) + \sum_{x \in \partial_- V_N} r(x) v(x) w(x),$$

where $r(x) := |\{y \in V_N^c : \operatorname{dist}(x, y) = 1\}|$. Obviously, $1 \leq r(x) \leq d$ for all $x \in \partial_- V_N$. It is immediate that $\mathcal{D}(\cdot, \cdot)$ is symmetric, bilinear and positive definite. We write $\|v\|_{\mathcal{D}} := \sqrt{\mathcal{D}(v, v)}$. In Appendix A we prove some estimates for discrete Sobolev norms and the Dirichlet form $\mathcal{D}(\cdot, \cdot)$.

To compare $G_N$ and $\overline{G}_N$, we use the fact that the difference of the two Green's functions,

$$H_N(x, y) := \overline{G}_N(x, y) - G_N(x, y),$$

satisfies the following boundary value problem:

$$\Delta^2 H_N(x, y) = 0, \qquad y \in V_N,$$
$$H_N(x, y) = \overline{G}_N(x, y), \qquad y \in \partial_2 V_N.$$

Let $f$ be any function $V_N \cup \partial_2 V_N \to \mathbb{R}$ which satisfies $f(y) = \overline{G}_N(x, y)$ for all $y \in \partial_2 V_N$. Then $u(y) := H_N(x, \cdot) - f(\cdot)$ satisfies

$$(5) \qquad \begin{aligned} \Delta^2 u(y) &= g(y), \qquad y \in V_N, \\ u(y) &= 0, \qquad y \in \partial_2 V_N, \end{aligned}$$

where $g(y) := -\Delta^2 f(y)$. The idea is now to choose an $f$ sufficiently regular in the interior of $V_N$, and show that this yields a solution $u$ of (5) which is $C^1$ in the discrete sense on $V_N^\delta$, meaning that if $x \in V_N^\delta, 0 < \delta < 1/2$, we have $\sup_{y \in V_N^\delta} |u(y)| \leq c$ and $\sup_{y \in V_N^\delta} |\nabla u(y)| \leq \frac{c}{N}$. Then we can derive estimates on $H_N(x, y)$ for $x, y \in V_N^\delta$.

Note that a function $u$ is a solution of (5) if and only if for any function $v : V_N \cup \partial_2 V_N \to \mathbb{R}$ it satisfies

$$\sum_{x \in V_N} \Delta^2 u(x) v(x) = \sum_{x \in V_N} g(x) v(x).$$

(Take $v = 1_x, x \in V_N$.) Summation by parts now shows that, since $u \in E_1$,

$$\sum_{x \in V_N} \Delta^2 u(x) v(x) = \mathcal{D}(u, v).$$

Hence, $\mathcal{D}(\cdot, \cdot)$ is the Dirichlet form corresponding to our boundary value problem and, therefore, an equivalent formulation of (5) is

$$(6) \qquad \mathcal{D}(u, v) = \langle g, v \rangle_{L_2(V_N)} \qquad \forall v \in E_1,$$

where $\langle \cdot, \cdot \rangle_{L_2(V_N)}$ denotes the $L_2$ scalar product on $V_N$. The Riesz Theorem now gives us a "weak" solution of (6): Clearly, for fixed $w \in E_1$, the map $v \mapsto \mathcal{D}(v, w)$ is well defined and linear from $E_1 \to \mathbb{R}$, so that by Riesz there exists $h_w \in E_1$ such that $\mathcal{D}(v, w) = \langle h_w, v \rangle_{L_2(V_N)}$, and the map $A: w \mapsto h_w$ is well defined and linear. It is injective, and therefore bijective since $E_1$ is finite dimensional. Thus, $A^{-1}$ exists, and $u := A^{-1}(-\Delta^2 f)$ is a solution of (6) and therefore also a solution of (5).

LEMMA 2.3.    *The unique solution $u$ of (5) satisfies*

$$\|u\|_{H^2(V_N)} \leq cN^4 \|g\|_{L_2(V_N)}.$$

PROOF.    We have just shown existence and uniqueness. For the norm estimate, note that by Corollary A.6 we have $\|u\|_{H^2(V_N)}^2 \leq cN^4 \mathcal{D}(u, u) = cN^4 \langle g, u \rangle_{L^2(V_N)} \leq cN^4 \|g\|_{L^2(V_N)} \|u\|_{L^2(V_N)}$. But this implies $\|u\|_{H^2(V_N)} \leq cN^4 \|g\|_{L^2(V_N)}$.    □

Let us now return to the case where $g = -\Delta^2 f$, where we want $f$ to satisfy the following:

LEMMA 2.4.    *Let $d = 4$. Let $0 < \delta < 1/2$, and $0 < \delta' < \delta/2$, and let $x \in V_N^\delta$. There exists a function $f$ on $V_N$ which satisfies the following conditions: There is a constant $c = c(\delta, \delta') > 0$ such that:*

(a)  $f(y) = \overline{G}_N(x, y)$ *for all $y \in V_N \setminus V_N^{\delta'}$,*

(b)  $|\nabla^\alpha f(y)| \leq \frac{c}{N^{|\alpha|}}$ *for all $y$ in $V_N^\delta$ and $|\alpha| \leq 5$,*

(c)  $|\Delta^2 f(y)| \leq \frac{c}{N^4}$, *and $|\nabla^i \Delta^2 f(y)| \leq \frac{c}{N^5}$ for all $y \in V_N$.*

PROOF.    It suffices to show that $|\nabla^\alpha \overline{G}_N(y)| \leq \frac{c}{N^{|\alpha|}}$ for all $y$ with $\delta' N \leq \mathrm{dist}(y, V_N^c) \leq (\delta/2)N$ and $|\alpha| \leq 5$. Then we can choose $f$ equal to any regular function on $V_N^\delta$, equal to $\overline{G}_N$ on $V_N \setminus V_N^{\delta'}$, and interpolate in between, which is possible since the number of interpolation points is of order $N^4$.

If $\alpha = (\alpha_1, \ldots, \alpha_d) \in \mathbb{N}_0^d$ and $f: \mathbb{R}^d \to \mathbb{R}$, we write $D^\alpha f(y) := \frac{\partial^{\alpha_1} \ldots \partial^{\alpha_d}}{\partial y_1^{\alpha_1} \ldots \partial y_d^{\alpha_d}} f(y)$. Note that the proof of Theorem 1.5.5 of [11] can be generalized to show that, if $y \neq 0$,

$$\nabla^\alpha \Gamma(0, y) = a_d D^\alpha(|y|^{2-d}) + O(|y|^{-d-|\alpha|+1})$$

for some constant $a_d$. Since $\Gamma_N(x,y) = \Gamma(x,y) - \sum_{z \in \partial V_N} \mathbb{P}^0(X_{\tau_N} = z)\Gamma(z,y)$, it follows immediately that for any $y$ with $\text{dist}(y, \partial V_N) \geq \delta' N$ and $|x - y| \geq (\delta/2)N$ we have

$$|\nabla^\alpha \Gamma_N(x,y)| \leq c(\delta, \delta')N^{-d-|\alpha|+2}.$$

We first assume $x = 0$. Split

$$\nabla^\alpha \overline{G}_N(0,y) = \sum_{z \in V_N} \Gamma_N(0,z)\nabla^\alpha \Gamma_N(z,y)$$

$$= \sum_{z \in V_N^\delta} \Gamma_N(0,z)\nabla^\alpha \Gamma_N(z,y) + \sum_{z \in V_N \setminus V_N^\delta} \Gamma_N(0,z)\nabla^\alpha \Gamma_N(z,y).$$

If $z \in V_N^\delta$ and $\text{dist}(y, V_N^c) \geq \delta' N$, we have $|z - y| \geq \delta' N$, and we can bound the first term by

$$\left| \sum_{z \in V_N^\delta} \Gamma_N(0,z)\nabla^\alpha \Gamma_N(z,y) \right| \leq \frac{c}{N^{d+|\alpha|-2}} \sum_{z \in V_N^\delta} \frac{1}{|z|^{d-2}} \leq \frac{c}{N^{d+|\alpha|-4}}.$$

The second term we split again:

$$\sum_{z \in V_N \setminus V_N^\delta} \Gamma_N(0,z)\nabla^\alpha \Gamma_N(z,y)$$

$$= \sum_{z \in V_N \setminus V_N^\delta} \Gamma_N(0,z)\nabla^\alpha \Gamma(z,y)$$

$$- \sum_{z \in V_N \setminus V_N^\delta} \sum_{w \in \partial V_N} \mathbb{P}^z(X_{\tau_N} = w)\Gamma_N(0,z)\nabla^\alpha \Gamma(w,y).$$

Again we have for any $w \in \partial V_N$ that $|w - y| \geq \delta' N$ and, therefore, as above,

$$\left| \sum_{z \in V_N \setminus V_N^\delta} \sum_{w \in \partial V_N} \mathbb{P}^z(X_{\tau_N} = w)\Gamma_N(0,z)\nabla^\alpha \Gamma(w,y) \right| \leq cN^{-d-|\alpha|+4}.$$

For the remaining term we use summation by parts (for $|\alpha| \leq 2$ this is not necessary, we could use similar estimates as before). Note that, since $\Gamma(z,y) = \Gamma(y,z)$, we have

$$\Gamma(z, y + e_i) - \Gamma(z,y) = \Gamma(z - e_i, y) - \Gamma(z,y)$$

and, thus,

$$\nabla^\alpha \Gamma(z,y) = \nabla^{-\alpha}\Gamma(y,z)$$

(we always let the difference operator act on the second variable). Thus, if $\alpha = \alpha' + e_i$, by summation by parts,

$$\sum_{z \in V_N \setminus V_N^\delta} \Gamma_N(0, z) \nabla^\alpha \Gamma(z, y)$$

$$= \sum_{z \in V_N \setminus V_N^\delta} \nabla^{e_i} \Gamma_N(0, z) \nabla^{\alpha'} \Gamma(z, y) + \sum_{z \in \partial(V_N \setminus V_N^\delta)} r(z) \Gamma_N(0, z) \nabla^{\alpha'} \Gamma(z, y),$$

where $1 \le r(z) \le d$ is the number of points in $V_N \setminus V_N^\delta$ which are neighbors of $z$. Note that

$$\sum_{z \in \partial(V_N \setminus V_N^\delta)} r(z) \Gamma_N(0, z) \nabla^{\alpha'} \Gamma(z, y) \le c N^{d-1} \frac{1}{N^{d-2}} \frac{1}{N^{d+|\alpha'|-2}} \le c \frac{1}{N^{d+|\alpha|-4}}.$$

Similarly, we have for any $\alpha', \beta$ with $|\alpha'| + |\beta| = |\alpha| - 1$ that

$$\sum_{z \in \partial(V_N \setminus V_N^\delta)} r(z) \nabla^\beta \Gamma_N(0, z) \nabla^{\alpha'} \Gamma(z, y) \le c \frac{1}{N^{d+|\alpha|-4}}.$$

Hence, we can iterate summation by parts and obtain that

$$\left| \sum_{z \in V_N \setminus V_N^\delta} \Gamma_N(0, z) \nabla^\alpha \Gamma(z, y) \right|$$

$$\le \left| \sum_{z \in V_N \setminus V_N^\delta} \nabla^\alpha \Gamma_N(0, z) \Gamma(z, y) \right| + c \frac{1}{N^{d+|\alpha|-4}}$$

$$\le c \frac{1}{N^{d+|\alpha|-2}} \sum_{z \in V_N \setminus V_N^\delta} \frac{1}{|z - y|^{d-2}} + c \frac{1}{N^{d+|\alpha|-4}}$$

$$\le \frac{c}{N^{d+|\alpha|-4}}.$$

This completes the proof, since similar arguments hold if $x \in V_N^\delta$ is arbitrary. $\square$

If we choose $f$ as in Lemma 2.4, we know from Lemma 2.3 that the solution $u$ of (5) is in $H^2(V_N)$ in the discrete sense:

COROLLARY 2.5.   *If* $\sup_{x \in V_N} |\Delta^2 f(x)| \le \frac{c}{N^4}$, *then* $\|u\|_{H^2(V_N)} \le N^{d/2}$.

For our purpose, we need stronger regularity of the solution than what we obtain from Lemma 2.3. To obtain this, we use a discrete version of the well-known bootstrap-technique in PDE; compare, for example, [15]. The first step is the following lemma.

LEMMA 2.6. *Let $1/2 < \delta < 1$, $0 < \varepsilon < 1/8$, and let $N$ be large enough, such that $\varepsilon N > 1$. Let $\chi \colon \mathbb{Z}^d \to \mathbb{R}$ satisfy $|\nabla^\alpha \chi| \leq cN^{-|\alpha|}$ for any multiindex $\alpha$, $\chi = 1$ on $V_N^\delta$ and $\chi(x) = 0$ if $\mathrm{dist}(x, \partial V_N) \leq 2\varepsilon N$. Furthermore, let $v \colon V_N \to \mathbb{R}$ be any function with $v(x) = 0$ if $\mathrm{dist}(x, \partial V_N) \leq \varepsilon N$. Then there exists $\overline{v}$ with $\|\overline{v}\|_{H^2(V_N)} = \|v\|_{H^2(V_N)}$, such that*

$$N^4 \mathcal{D}(N\nabla_i(\chi u), v) = -N^4 \langle g, N\chi \nabla_i \overline{v}\rangle_{L_2(V_N)} + I_0,$$

*where $I_0 \leq c\|u\|_{H^2(V_N)}\|v\|_{H^2(V_N)}$.*

PROOF. First, note the product rule for $\nabla_i$: $\nabla_i(vw)(x) = \nabla_i v(x)w(x) + v(x + e_i)\nabla_i w(x)$. Furthermore, if $v$ has support in the interior of $V_N$, then $\sum_{x \in V_N} \nabla_i v(x) = 0$. Using this and the assumptions on $v$, we get

$$\begin{aligned}
N^4 \mathcal{D}(N\nabla_i(\chi u), v) &= N^4 \sum_{x \in V_N} \Delta N\nabla_i(\chi u)(x)\Delta v(x) \\
&= N^4 \sum_{x \in V_N} N\nabla_i \Delta(\chi u)(x)\Delta v(x) \\
&= N^4 \sum_{x \in V_N} N\nabla_i(\Delta(\chi u)\Delta v)(x) \\
&\quad - N^4 \sum_{x \in V_N} (\Delta(\chi u))(x + e_i)N\nabla_i\Delta v(x).
\end{aligned}$$

Now the first term is $0$ due to the choice of the support of $v$, and the second—using the product rule on the discrete Laplacian—is equal to

$$-N^4 \sum_{x \in V_N} \Delta u(x + e_i)\chi(x + e_i)N\nabla_i\Delta v(x)$$

$$+ N^4 \sum_{x \in V_N} \sum_{\alpha \colon |\alpha| \leq 2} \sum_{\substack{\beta \colon |\beta| \leq 1 \\ |\alpha| + |\beta| = 2}} k(\alpha, \beta)(\nabla^\alpha \chi)(x + e_i)(\nabla^\beta u)(x + e_i)N\nabla_i\Delta v(x)$$

for suitable $k(\alpha, \beta) \in \mathbb{R}$. In the second term we use summation by parts and the regularity of $\chi$ to bound its absolute value by $c\|u\|_{H^2(V_N)}\|v\|_{H^2(V_N)}$. If we define the translation operator $\tau_i$ by $\tau_i(x) := x + e_i$, we can again use the product rule to rewrite the first term as

$$-N^4 \sum_{x \in V_N} \Delta u(x + e_i)\chi(x + e_i)N\nabla_i\Delta v(x)$$

$$= -N^4 \sum_{x \in V_N} (\Delta u)(x + e_i)\Delta((\chi \circ \tau_i)N\nabla_i v)(x)$$

$$+ \sum_{x \in V_N} (\Delta u)(x + e_i) \sum_{\alpha \colon |\alpha| \leq 2} \sum_{\substack{\beta \colon |\beta| \leq 1 \\ |\alpha| + |\beta| = 2}} k(\alpha, \beta)\nabla^\alpha \chi(x)\nabla^\beta N\nabla_i v(x).$$

Here, by (6), the first term is equal to

$$-N^4 \mathcal{D}(u, \chi N \nabla_i (v \circ \tau^{-1})) = -N^4 \langle g, \chi N \nabla_i (v \circ \tau^{-1}) \rangle_{L_2(V_N)},$$

and the second is again bounded from above by $c\|u\|_{H^2(V_N)} \|v\|_{H^2(V_N)}$.  □

PROPOSITION 2.7.   *Let $\chi$ as in Lemma 2.6, and let $u$ be the solution of (5) where $f$ satisfies the properties* (a), (b) *and* (c) *of Lemma 2.4. Then there exists $c > 0$ such that*

$$\|\chi u\|_{H^3(V_N)} \leq cN^{d/2}.$$

PROOF.   Let $\overline{v}$ be the same as in Lemma 2.6. Note that

$$|\langle g, N\chi \nabla_i \overline{v} \rangle_{L_2(V_N)}| \leq \|g\|_{L_2(V_N)} \|N\chi \nabla_i \overline{v}\|_{L_2(V_N)} \leq c\|g\|_{L_2(V_N)} \|\overline{v}\|_{H^1(V_N)}$$

$$\leq c\|g\|_{L_2(V_N)} \|v\|_{H^2(V_N)}.$$

Thus, if we set $v = N\nabla_i(\chi u)$ in Lemma 2.6 , we have, using Corollary A.6,

$$\|N\nabla_i(\chi u)\|_{H^2(V_N)}^2 \leq c_1 N^4 \mathcal{D}(N\nabla_i(\chi u), N\nabla_i(\chi u))$$

$$\leq c_1 \|N\nabla_i(\chi u)\|_{H^2(V_N)} (N^4 \|g\|_{L_2(V_N)} + \|u\|_{H^2(V_N)})$$

and so

$$\|N\nabla_i(\chi u)\|_{H^2(V_N)} \leq c(N^4 \|g\|_{L_2(V_N)} + \|u\|_{H^2(V_N)}) \leq cN^{d/2}$$

by Corollary A.6 and Lemma 2.3. The claim now follows from Remark A.5.  □

COROLLARY 2.8.   *Let $d = 4$. If $u$ is a solution of (5), where $f$ satisfies the properties* (a), (b) *and* (c) *of Lemma 2.4, and $\chi$ is defined as in Lemma 2.6, then $\chi u \in H^k(V_N)$ for $0 \leq k \leq 4$.*

PROOF.   Apply the arguments of Lemma 2.6 and Proposition 2.7 with $N\nabla_i u$ in the place of $u$, and $N\nabla_i g$ in the place of $g$, and use the result of Proposition 2.7.  □

Now we can conclude:

COROLLARY 2.9.   *Let $d = 4$, and $0 < \delta < 1/2$. There exists $c(\delta) > 0$ such that, for all $x \in V_N^\delta$,*

$$\sup_{y \in V_N^\delta} |\overline{G}_N(x, y) - G_N(x, y)| \leq c(\delta)$$

*and, for all $1 \leq i \leq d$,*

$$\sup_{y \in V_N^\delta} |\nabla_i(\overline{G}_N(x, y) - G_N(x, y))| \leq c(\delta)N^{-1}.$$

PROOF. By Corollary 2.8, $\chi u \in H^4(V_N)$ and, thus, by Corollary B.2, $\sup |\chi u| \le c$ and $\sup |\nabla_i \chi u| \le c/N$. Since $\chi = 1$ and $\nabla_i \chi = 0$ on $V_N^\delta$, this implies $\sup_{x \in V_N^\delta} |u(x)| \le c$ and $\sup_{x \in V_N^\delta} |\nabla_i u(x)| \le c/N$. Since $\overline{G}_N(x, y) - G_N(x, y) = u(y) + f(y)$, the claim is proven by the assumptions we made on $f$. □

Corollary 2.9, together with Lemma 2.2, finally proves the logarithmic variance structure of the membrane model, which proves Proposition 1.1.

PROOF OF PROPOSITION 1.1. Note that $\mathrm{var}_N(\varphi_x) \le \mathrm{var}_N(\varphi_0)$ for all $x \in V_N$. Then both claims follow from the estimates on $\overline{G}_N$ in Lemma 2.2 and Corollary 2.9. □

Additionally to Proposition 1.1, Lemma 2.11 below will be crucial for the approximation of the field with a hierarchical one (see [1]). We therefore introduce the discrete version of the fundamental solution for the Bilaplacian: Let, as before, $(X_k)_{k \in \mathbb{N}}$ be a simple random walk on the lattice, and let $\mathbb{P}^x$ denote it's law conditional on starting in $x$. Let

$$a(x, y) := \sum_{k=0}^\infty (k+1)(\mathbb{P}^x(X_k = x) - \mathbb{P}^x(X_k = y)).$$

Lemma 2.10 below shows that this is finite for any pair $x, y \in \mathbb{Z}^d$. Note first that $a(0, 0) = 0$, and that $a(x, y) = a(0, y - x)$. The local central limit theorem ([11], Theorem 1.2.1) allows us to compute $a(x, y)$:

LEMMA 2.10. Let $d = 4$. There exists a constant $K$, such that for all $y \ne 0$, for all $0 < \alpha < 2$,

$$(7) \qquad a(0, y) = \frac{8}{\pi^2} \log |y| + K + o(|y|^{-\alpha}).$$

PROOF. First, note that $a(0, y) = \sum_{k=0}^\infty k(\mathbb{P}^0(X_k = 0) - \mathbb{P}^0(X_k = y)) + \Gamma(0, 0) - \Gamma(0, y)$. Remember that $\Gamma(0, y) \le O(|y|^{-2})$, and $\Gamma(0, 0)$ is a constant. Let $\overline{p}(k, x) := \frac{8}{\pi^2 k^2} \exp(-\frac{2|x|^2}{k})$ and

$$E(k, x) := \begin{cases} \mathbb{P}^0(X_k = x) - \overline{p}(k, x), & \text{if } \mathbb{P}^0(X_k = x) \ne 0, \\ 0, & \text{otherwise.} \end{cases}$$

Let us first assume that $y$ is even. Then

$$\sum_{k=0}^\infty k(\mathbb{P}^0(X_k = 0) - \mathbb{P}^0(X_k = y)) = \sum_{k=1}^\infty 2k(\mathbb{P}^0(X_{2k} = 0) - \mathbb{P}^0(X_{2k} = y))$$

and

$$\sum_{k=1}^{\infty} 2k(\mathbb{P}^0(X_{2k}=0) - \mathbb{P}^0(X_{2k}=y))$$

$$= \sum_{k=1}^{\infty} 2k(\overline{p}(2k,0) - \overline{p}(2k,y) + E(2k,0) - E(2k,y)).$$

We first consider the remainder term. From the local CLT with error bounds ([11], Theorem 1.2.1) we know

$$|E(k,y)| \le O(k^{-3}) \quad \text{and} \quad |E(k,y)| \le |y|^{-2}O(k^{-2})$$

and, consequently,

$$\sum_{k=1}^{\infty} 2kE(2k,y) \le \sum_{k \le |y|^2/2} 2kE(2k,y) + \sum_{k > |y|^2/2} 2kE(2k,y)$$

$$\le |y|^2 \sum_{k \le |y|^2/2} E(2k,y) + \sum_{k > |y|^2/2} 2kO((2k)^{-3})$$

$$\le |y|^2 \sum_{k \le |y|^2/2} E(2k,y) + O(|y|^{-2}).$$

But from Lemma 1.5.2 of [11] we know that $\sum_{k=0}^{\infty} E(k,y) = o(|y|^{-\alpha})$ for any $\alpha < 4$ as $|y| \to \infty$.

Now consider the other term. By definition,

$$\sum_{k=1}^{\infty} 2k(\overline{p}(2k,0) - \overline{p}(2k,y)) = \frac{4}{\pi^2} \sum_{k=1}^{\infty} \frac{1}{k}(1 - \exp(-|y|^2/k)).$$

Now use exactly the same steps as in the proof of Theorem 1.6.2 of [11] to show that there is a constant $\tilde{K}$ such that

$$\frac{4}{\pi^2} \sum_{k=1}^{\infty} \frac{1}{k}(1 - \exp(-|y|^2/k)) = \frac{4}{\pi^2}(\log|y|^2 + \tilde{K} + O(|y|^{-2})).$$

This proves the case where $y$ is even with $K = \Gamma(0,0) + \frac{4}{\pi^2}\tilde{K} + \sum_{k=1}^{\infty} 2kE(2k,0)$. If $y$ is odd,

$$\sum_{k=0}^{\infty} k(\mathbb{P}^0(X_k=0) - \mathbb{P}^0(X_k=y))$$

$$= \sum_{k=1}^{\infty} 2k(\mathbb{P}^0(X_{2k}=0) - \mathbb{P}^0(X_{2k+1}=y)) - \Gamma(0,y)$$

$$= \frac{1}{2d} \sum_{v:|y-v|=1} \sum_{k=1}^{\infty} 2k(\mathbb{P}^0(X_{2k}=0) - \mathbb{P}^0(X_{2k}=v)) - \Gamma(0,y).$$

Of course, all these $v$ are even, so we obtain, since $\frac{1}{2d}\sum_{v:|y-v|=1}\log|v|^2 = \log|y|^2 + O(|y|^{-2})$,

$$a(0,y) = \frac{4}{\pi^2}\frac{1}{2d}\sum_{v:|y-v|=1}\log|v|^2 + K + o(|y|^{-\alpha}) = \frac{8}{\pi^2}\log|y| + K + O(|y|^{-\alpha}),$$

where $\alpha < 2$ and $K$ is the same as before. $\quad\square$

This result together with the random walk representation for $\overline{G}_N$ is the key to proving the following result:

LEMMA 2.11. *Letting $0 < n < N$, let $A_N \subset \mathbb{Z}^d$ be a box of side-length $N$ and $A_n \subset A_N$ be a box of side-length $n$ with the same center $x_B \in \mathbb{Z}^d$ as $A_N$. Let $0 < \varepsilon < 1/2$. There exists $c > 0$ such that, for all $x \in A_n$ with $|x - x_B| \le \varepsilon n$,*

$$\mathrm{var}(E(\varphi_x|\mathcal{F}_{\partial_2 A_n}) - E(\varphi_{x_B}|\mathcal{F}_{\partial_2 A_n})|\mathcal{F}_{\partial_2 A_N}) \le c\varepsilon.$$

PROOF. Note that for any two subsets $E \subset F$ of $\mathbb{Z}^d$ we have

$$(8) \quad \mathrm{var}(\varphi_x|\mathcal{F}_{F^c}) = \mathrm{var}(\varphi_x|\mathcal{F}_{E^c}) + \mathrm{var}(E(\varphi_x|\mathcal{F}_{E^c})|\mathcal{F}_{F^c}) \ge \mathrm{var}(\varphi_x|\mathcal{F}_{E^c}).$$

Let $B_n := B_n(x_B) = \{z \in \mathbb{Z}^d : |x_B - z| < n\}$ be the ball of radius $n$ around $x_B$. We define $G_{B_n}$ analogous to $G_N$ as the Green's function of the biharmonic problem (1) on $B_n$ instead of $V_N$. Likewise, $\overline{G}_{B_n}$ is defined by (3) and (4) on $B_n$, and $H_{B_n} := G_{B_n} - \overline{G}_{B_n}$. It is clear that the regularity considerations of this section apply to $G_{B_n}$ and $\overline{G}_{B_n}$ as well and, thus, Corollary 2.9 can be applied. Note $B_n \subset A_n$, and so

$$\begin{aligned}
&\mathrm{var}(E(\varphi_x - \varphi_{x_B}|\mathcal{F}_{\partial_2 A_n})|\mathcal{F}_{\partial_2 A_N}) \\
&= \mathrm{var}(\varphi_x - \varphi_{x_B}|\mathcal{F}_{A_N^c}) - \mathrm{var}(\varphi_x - \varphi_{x_B}|\mathcal{F}_{A_n^c}) \\
(9) \quad &\le \lim_{N\to\infty}(\mathrm{var}(\varphi_x - \varphi_{x_B}|\mathcal{F}_{A_N^c}) - \mathrm{var}(\varphi_x - \varphi_{x_B}|\mathcal{F}_{B_n^c})) \\
&= \lim_{N\to\infty}(G_N(x,x) - 2G_N(x,x_B) + G_N(x_B,x_B) \\
&\quad - G_{B_n}(x,x) + 2G_{B_n}(x,x_B) - G_{B_n}(x_B,x_B)).
\end{aligned}$$

(Of course we do not know if the limit exists, but otherwise the rhs is equal to $+\infty$.) Now, $G_N = \overline{G}_N + H_N$. From Corollary 2.9 we know that $|H_N(y,z) - H_N(y, z + e_i)| \le cN^{-1}$, and since $|x - x_B| \le \varepsilon n$, we need at most $4\varepsilon n$ steps to get from $x_B$ to $x$. Thus, $|H_N(y,x) - H_N(y,x_B)| \le \varepsilon n \cdot cN^{-1}$ if $y \in \{x, x_B\}$, and so

$$\begin{aligned}
&\lim_{N\to\infty}(H_N(x,x) - 2H_N(x,x_B) + H_N(x_B,x_B) \\
&\quad - H_{B_n}(x,x) + 2H_{B_n}(x,x_B) - H_{B_n}(x_B,x_B)) \\
&\quad \le \lim_{N\to\infty}\varepsilon n \cdot cN^{-1} + \varepsilon n \cdot cn^{-1} \le c\varepsilon.
\end{aligned}$$

We are therefore left with estimating the terms in (9) involving $\overline{G}_N$ and $\overline{G}_{B_n}$. We have

$$\overline{G}_N(x,x) - 2\overline{G}_N(x,x_B) + \overline{G}_N(x_B,x_B)$$
$$\qquad - \overline{G}_{B_n}(x,x) + 2\overline{G}_{B_n}(x,x_B) - \overline{G}_{B_n}(x_B,x_B)$$
$$= \sum_{k=0}^{\infty}(k+1)[\mathbb{P}^x(X_k = x, \tau_{B_n} \leq k \leq \tau_{B_N}) - \mathbb{P}^x(X_k = x_B, \tau_{B_n} \leq k \leq \tau_{B_N})$$
$$\qquad\qquad + \mathbb{P}^{x_B}(X_k = x_B, \tau_{B_n} \leq k \leq \tau_{B_N})$$
$$\qquad\qquad\qquad\qquad - \mathbb{P}^{x_B}(X_k = x, \tau_{B_n} \leq k \leq \tau_{B_N})].$$

Hence, using the above monotonicity (8), we are done if we show

$$(10) \quad \sum_{k=0}^{\infty}(k+1)[\mathbb{P}^x(X_k = x, k \geq \tau_{B_n}) - \mathbb{P}^x(X_k = x_B, k \geq \tau_{B_n})$$
$$\qquad\qquad + \mathbb{P}^{x_B}(X_k = x_B, k \geq \tau_{B_n}) - \mathbb{P}^{x_B}(X_k = x, k \geq \tau_{B_n})] \leq c\varepsilon.$$

Define

$$T_1 := \sum_{z \in \partial B_n}(\mathbb{P}^x(X_{\tau_{B_n}} = z) - \mathbb{P}^{x_B}(X_{\tau_{B_n}} = z))(a(z,x_B) - a(z,x))$$

and

$$T_2 := \sum_{z \in \partial B_n}\sum_{m=0}^{\infty}m(\mathbb{P}^x(\tau_{B_n} = m, X_{\tau_{B_n}} = z) - \mathbb{P}^{x_B}(\tau_{B_n} = m, X_{\tau_{B_n}} = z))$$
$$\qquad\qquad \times (\Gamma(z,x) - \Gamma(z,x_B)).$$

Due to Lemma 2.10, for $x, x_B$ as above, $\sup_{z \in \partial B_n}|a(z,x_B) - a(z,x)| \leq c\varepsilon$, which implies $|T_1| \leq c\varepsilon$. For $T_2$, observe that, by construction, $|z - x_B| \geq n$ and $|z - x| \geq (1-\varepsilon)n$, which implies $\sup_{z \in \partial B_n}\Gamma(z,x) \leq \frac{c}{(1-\varepsilon)^2 n^2}$ and likewise for $\Gamma(z,x_B)$. On the other hand,

$$\sum_{z \in \partial B_n}\sum_{m=0}^{\infty}m(\mathbb{P}^x(\tau_{B_n} = m, X_{\tau_{B_n}} = z) - \mathbb{P}^{x_B}(\tau_{B_n} = m, X_{\tau_{B_n}} = z))$$
$$\qquad\qquad = \mathbb{E}^x(\tau_{B_n}) - \mathbb{E}^{x_B}(\tau_{B_n}).$$

From [11], Equation 1.21, we know that

$$n^2 - |y - x_B|^2 \leq \mathbb{E}^y(\tau_{B_n}) \leq (n+1)^2 - |y - x_B|^2$$

for all $y \in B_n$. Therefore, $|\mathbb{E}^x(\tau_{B_n}) - \mathbb{E}^{x_B}(\tau_{B_n})| \leq \varepsilon^2 n^2 + 2n + 1$, and if $n$ is large enough, $|T_2| \leq c\varepsilon$. Thus, we have shown

$$(11) \qquad\qquad\qquad |T_1 + T_2| \leq c\varepsilon$$

for some finite $c$. We have by definition of $\Gamma(\cdot, \cdot)$ and $a(\cdot, \cdot)$,

$$T_1 + T_2 = \sum_{k=0}^{\infty} \sum_{m=0}^{\infty} \sum_{z \in \partial B_n} (k+m+1)(\mathbb{P}^z(X_k = x) - \mathbb{P}^z(X_k = x_B)) \tag{12}$$

$$\times (\mathbb{P}^x(\tau_{B_n} = m, X_{\tau_{B_n}} = z) - \mathbb{P}^{x_B}(\tau_{B_n} = m, X_{\tau_{B_n}} = z)).$$

By the Markov property,

$$\mathbb{P}^x(X_k = x, k \geq \tau_{B_n}) \tag{13}$$

$$= \sum_{m=0}^{\infty} \sum_{z \in \partial B_n} \mathbb{P}^z(X_{k-m} = x)\mathbb{P}^x(\tau_{B_n} = m, X_{\tau_{B_n}} = z)$$

and similarly for $\mathbb{P}^x(X_k = x_B, k \geq \tau_{B_n})$ etc. Equations (12) and (13) imply

$$\sum_{k=0}^{\infty}(k+1)[\mathbb{P}^x(X_k = x, k \geq \tau_{B_n}) - \mathbb{P}^x(X_k = x_B, k \geq \tau_{B_n})$$

$$+ \mathbb{P}^{x_B}(X_k = x_B, k \geq \tau_{B_n}) - \mathbb{P}^{x_B}(X_k = x, k \geq \tau_{B_n})] \leq T_1 + T_2 \leq c\varepsilon,$$

the last inequality by (11). This completes the proof of (10). □

**3. Maximum of the field.** In this section we prove Theorem 1.2, using the strategy of [1] and [5], whose crucial ingredients are the logarithmic structure of the variances (Proposition 1.1) and the concentration result (Lemma 2.11). Let $\alpha \in (1/2, 1)$. We cover $V_N^\delta$ with boxes of side-length $N^\alpha$ as in [1]: Let $x_0 \in V_N$, and let

$$M_\alpha := \{x_0 + i(N^\alpha + 2) : i = (i_1, \ldots, i_4) \in \mathbb{N}^4 \text{ such that } x_0 + i(N^\alpha + 2) \subset V_N\}.$$

We consider the set of boxes $B$ with midpoint in $M_\alpha$ and side-length $N^\alpha$. We will always assume that $N^\alpha$ is an odd integer, which is no restriction as $N \to \infty$. By construction, the boundaries between two boxes have thickness 2 (on the lattice), which is the range of interactions of $\Delta^2$. Let $\Pi_\alpha$ denote the set of such boxes which are contained in $V_N^\delta$, and let $\Lambda_\alpha := \bigcup_{B \in \Pi_\alpha} \partial_2 B$ be the set of all boundaries of boxes in $\Pi_\alpha$. We denote by $\mathcal{F}_\alpha$ the sigma-algebra generated by the $\varphi_x : x \in \Lambda_\alpha$. Conditional on $\mathcal{F}_\alpha$, what happens inside different boxes is independent.

Now fix $K \in \mathbb{N}$. Set $\alpha_i := \alpha(1 - \frac{i-1}{K}), 1 \leq i \leq K+1$. We define the following sets of boxes: First, let $\Gamma_{\alpha_1} := \Pi_{\alpha_1}$. Then $\Gamma_{\alpha_i}, i \geq 2$, is defined recursively: For $B \in \Gamma_{\alpha_{i-1}}$, let $\Gamma_{B, \alpha_i} := \{B' \in \Pi_{\alpha_i} : B' \subset B/2\}$, and $\Gamma_{\alpha_i} := \bigcup_{B \in \Gamma_{\alpha_{i-1}}} \Gamma_{B, \alpha_i}$. For $B \in \Pi_\alpha$, we denote the midpoint of $B$ by $x_B$. Let

$$\varphi_B := E_N(\varphi_{x_B} | \mathcal{F}_{\partial_2 B}) = E_N(\varphi_{x_B} | \mathcal{F}_\alpha).$$

If $B \in \Pi_{\alpha_i}$ and $B' \in \Pi_{\alpha_j}$, with $\alpha_i \leq \alpha_j$ such that $x_B = x_{B'}$, by (8) and Proposition 1.1, we see that

$$\text{var}(\varphi_B | \mathcal{F}_{\alpha_j}) = \text{var}(\varphi_{x_B} | \mathcal{F}_{\alpha_j}) - \text{var}(\varphi_{x_B} | \mathcal{F}_{\alpha_i}) = \gamma(\alpha_j - \alpha_i) \log N + O(1). \quad (14)$$

Note that, by (2), there exist coefficients $h(z) \in \mathbb{R}$ such that

$$\varphi_B = \sum_{z \in \partial_2 B} h(z) \varphi_z.$$

Unlike in the case of the lattice free field, however, the $h(z)$ need not lie between 0 and 1 (in fact, one can see that there are both positive and negative coefficients, and they need not be bounded in $N$). Some arguments in the proof need to be adapted to this fact, in particular, comparing $\varphi_B$ and $\varphi_{x_B}$ requires some work, for which we use Gaussian tail estimates. For the sake of readability, we give a complete proof, including also those parts that are practically identical to [1] or [5]. Note that one direction is easy to prove:

PROOF OF THEOREM 1.2(A). Using Proposition 1.1, we obtain

$$P_N \left( \sup_{x \in V_N} \varphi_x \geq 2\sqrt{2\gamma} \log N \right)$$

$$\leq |V_N| \sup_{x \in V_N} P_N(\varphi_x \geq 2\sqrt{2\gamma} \log N)$$

$$\leq N^4 \frac{\sqrt{\gamma \log N + c}}{\sqrt{2\pi} 2\sqrt{2\gamma} \log N} \exp\left( -\frac{(2\sqrt{2\gamma} \log N)^2}{2\gamma \log N + O(1)} \right),$$

which tends to zero as $N \to \infty$.   □

The second part is obtained from the following more general result (compare [5]):

THEOREM 3.1. *Let $0 < \delta < 1/2$, and let $0 < \lambda_0 < 1$ and $\lambda_0 < \lambda < 1$. For all $\varepsilon > 0$, there exists $c = c(\delta, \lambda_0) > 0$ such that*

$$P_N(|\{x \in V_N^\delta : \varphi_x \geq 2\sqrt{2\gamma}\lambda \log N\}| \leq N^{4(1-\lambda^2)-\varepsilon}) \leq \exp(-c(\log N)^2).$$

PROOF OF THEOREM 1.2(B). Chose in Theorem 3.1 $\lambda$ sufficiently close to 1, such that $2\sqrt{2\gamma}\lambda \geq (2\sqrt{2\gamma} - \eta)$ and $4\lambda^2 > 4 - \varepsilon$ are both satisfied.   □

To prove Theorem 3.1, we start on level $\alpha = \alpha_1$ of the box structure introduced before, and show that, on this level, a sufficiently high number of the $\varphi_B, B \in \Gamma_\alpha$, are positive.

LEMMA 3.2.    *Let $1/2 < \delta < 1$ and $\alpha \in (1/2, 1)$. There exist positive constants $\kappa, a$ depending on $\alpha$ and $\delta$, such that*

$$P_N(|\{B \in \Gamma_\alpha : \varphi_B \geq 0\}| \leq N^\kappa) \leq \exp(-a(\log N)^2).$$

PROOF.    Set $\alpha' = (1+\alpha)/2$, which implies $\alpha' > \alpha$. We consider the event

$$A := \left\{ \sharp \left\{ B \in \Pi_{\alpha'} : \varphi_B \geq \frac{-(1-\alpha')\sqrt{2\gamma}\log N}{2} \right\} \geq N^{1-\alpha'} \right\}.$$

The lemma will be proven showing that the following two estimates hold:

$$(15) \qquad P_N(A \cap \{\sharp\{B \in \Gamma_\alpha : \varphi_B \geq 0\} \leq N^\kappa\}) \leq \exp(-c(\log N)^2)$$

for some $c > 0$, and

$$(16) \qquad\qquad P_N(A^c) \leq \exp(-c(\log N)^2).$$

Obviously, these two estimates prove the lemma. We start with the second estimate. Let us split the event $A^c$ into

$$(17) \qquad \begin{aligned} P_N(A^c) &\leq P_N\left( A^c \cap \left\{ \max_{B \in \Pi_{\alpha'}} \varphi_B \leq (\log N)^2 \right\} \right) \\ &\quad + P_N\left( \max_{B \in \Pi_{\alpha'}} \varphi_B > (\log N)^2 \right) \end{aligned}$$

and bound the two terms. First, notice that for any $0 < \rho < 1$ we have

$$(18) \qquad \begin{aligned} P_N &\left( \max_{x \in V_N} \varphi_x > (1-\rho)(\log N)^2 \right) \\ &\leq N^4 \max_{x \in V_N} P_N(\varphi_x > (1-\rho)(\log N)^2) \\ &\leq N^4 \exp\left( -\frac{(1-\rho)^2(\log N)^4}{2\gamma \log N + C} \right) \\ &\leq \exp(-c(\log N)^3). \end{aligned}$$

Now we get

$$(19) \qquad \begin{aligned} P_N &\left( \left\{ \max_{B \in \Pi_{\alpha'}} \varphi_B > (\log N)^2 \right\} \cap \left\{ \max_{x \in V_N} \varphi_x \leq (1-\rho)(\log N)^2 \right\} \right) \\ &\leq P_N\left( \left\{ \max_{B \in \Pi_{\alpha'}} \varphi_B > (\log N)^2 \right\} \cap \left\{ \max_{x \in \Pi_{\alpha'}} \varphi_{x_B} \leq (1-\rho)(\log N)^2 \right\} \right) \\ &\leq |\Pi_{\alpha'}| \max_{B \in \Pi_{\alpha'}} P_N(\{\varphi_B > (\log N)^2\} \cap \{\varphi_{x_B} \leq (1-\rho)(\log N)^2\}) \\ &\leq cN^4 E_N(P_N(\varphi_{x_{B_0}} \leq (1-\rho)(\log N)^2 | \mathcal{F}_{\partial B_0}) 1_{\{\varphi_{B_0} > (\log N)^2\}}) \end{aligned}$$

for some fixed $B_0 \in \Pi_{\alpha'}$. Since by Proposition 1.1, conditional on $\mathcal{F}_{B_0}$, the random variable $\varphi_{x_{B_0}} - \varphi_{B_0}$ is centered Gaussian with $\mathrm{var}(\varphi_{x_{B_0}} - \varphi_{B_0}) \leq \gamma \alpha' \log N$, we have on $\{\varphi_{B_0} > (\log N)^2\}$

$$
\begin{aligned}
& P_N(\varphi_{x_{B_0}} \leq (1-\rho)(\log N)^2 | \mathcal{F}_{\partial B_0}) \\
(20) \quad & \leq P_N(\varphi_{x_{B_0}} - \varphi_{B_0} \leq -\rho(\log N)^2 | \mathcal{F}_{\partial B_0}) \\
& \leq \exp(-c(\log N)^3).
\end{aligned}
$$

Together, (18), (19) and (20) give the required bound on the second term in (17). To bound the first term, note that on $A^c \cap \{\max_{B \in \Pi_{\alpha'}} \varphi_B \leq (\log N)^2\}$ we have

$$
\begin{aligned}
& \frac{1}{|\Pi_{\alpha'}|} \sum_{B \in \Pi_{\alpha'}} \varphi_B \\
(21) \quad & \leq \frac{-(1-\alpha')\sqrt{2\gamma} \log N}{2} + \frac{N^{1-\alpha'}}{|\Pi_{\alpha'}|}\left(\frac{(1-\alpha')\sqrt{2\gamma} \log N}{2} + (\log N)^2\right).
\end{aligned}
$$

Since $|\Pi_{\alpha'}| = O(N^{4(1-\alpha')})$, we get from (21)

$$
(22) \quad \frac{1}{|\Pi_{\alpha'}|} \sum_{B \in \Pi_{\alpha'}} \varphi_B \leq \frac{-(1-\alpha')\sqrt{2\gamma} \log N}{3}.
$$

By Lemma C.1, we know that $\mathrm{var}(\frac{1}{|\Pi_{\alpha'}|} \sum_{B \in \Pi_{\alpha'}} \varphi_B) < \infty$, therefore, we obtain with (22)

$$
\begin{aligned}
& P_N\left(A^c \cap \left\{\max_{B \in \Pi_{\alpha'}} | \varphi_B \leq (\log N)^2\right\}\right) \\
& \leq P_N\left(\frac{1}{|\Pi_{\alpha'}|} \sum_{B \in \Pi_{\alpha'}} \varphi_B \leq \frac{-(1-\alpha')\sqrt{2\gamma} \log N}{3}\right) \\
& \leq \exp\left(\frac{-(1-\alpha')^2 \gamma (\log N)^2}{9 \mathrm{var}(1/|\Pi_{\alpha'}| \sum_{B \in \Pi_{\alpha'}} \varphi_B)}\right) \\
& \leq \exp(-c(\log N)^2).
\end{aligned}
$$

This gives the second bound in (17) and thus proves (16). For the proof of (15), we consider only the set of boxes in $\Pi_\alpha$ which have the same center as some box of $\Pi_{\alpha'}$: Let

$$
\Pi_{\alpha,\alpha'} := \{B \in \Pi_\alpha : \exists B' \in \Pi_{\alpha'} \text{ with } x_B = x_{B'}\}.
$$

We have

$$
P_N(A \cap \{|\{B \in \Gamma_\alpha : \varphi_B \geq 0\}| \leq N^\kappa\})
$$

$$(23) \qquad \leq P_N(A \cap \{|\{B \in \Pi_{\alpha,\alpha'} : \varphi_B \geq 0\}| \leq N^\kappa\})$$

$$\leq E_N(P_N(|\{B \in \Pi_{\alpha,\alpha'} : \varphi_B \geq 0\}| \leq N^\kappa | \mathcal{F}_{\alpha'}) 1_A).$$

We know that on $A$ there exist at least $N^{1-\alpha'}$ boxes $B' \in \Pi_{\alpha'}$ where there is $\varphi_{B'} \geq -(1-\alpha')\sqrt{2\gamma}\log N/2$. Choose $N^{1-\alpha'}$ of them and call them $B'_1, \dots,$ $B'_{N^{1-\alpha'}}$. Let $B_i \in \Pi_{\alpha,\alpha'}$ be the box with center $x_{B_i} = x_{B'_i}$. Set

$$\zeta_i := \varphi_{B_i} - \varphi_{B'_i}.$$

Then for $\kappa < 1 - \alpha'$,

$$(24) \qquad P_N(|\{B \in \Pi_{\alpha,\alpha'} : \varphi_B \geq 0\}| \leq N^\kappa | \mathcal{F}_{\alpha'})$$

$$\leq P_N\left( \sum_{i=1}^{N^{1-\alpha'}} 1_{\{\zeta_i \geq (1-\alpha')\sqrt{2\gamma}(\log N)/2\}} \leq N^\kappa \right).$$

By construction, we have $\varphi_{B'_i} = E_N(\varphi_{x_{B'_i}} | \mathcal{F}_{\alpha'}) = E_N(E_N(\varphi_{x_{B_i}} | \mathcal{F}_\alpha) | \mathcal{F}_{\alpha'}) = E_N(\varphi_{B_i} | \mathcal{F}'_\alpha)$. Therefore, we know the following:

- The $\zeta_i$ are centered Gaussian random variables under $P_N(\cdot | \mathcal{F}_{\alpha'})$.
- By (14), $\text{var}(\zeta_i) = \text{var}_{B'_i}(\varphi_{B_i}) = \gamma(1 - \alpha') \log N + O(1)$, since $\alpha' - \alpha = 1 - \alpha'$.

This implies

$$(25) \quad P_N\left( \zeta_i \geq \frac{1-\alpha'}{2}\sqrt{2\gamma}\log N \right) \geq \exp\left( \frac{-(1-\alpha')\log N}{4} \right) = N^{-(1-\alpha')/4}.$$

If we choose now $\kappa = (1 - \alpha')/2$ and set $\theta_i = 1_{\{\zeta_i \geq (1-\alpha')\sqrt{2\gamma}(\log N)/2\}}$, we know that on $A$ we have $\sum_{i=1}^{N^{1-\alpha'}} \theta_i \leq N^{(1-\alpha')/2}$, and from (25) we get $E\theta_i \geq N^{-(1-\alpha')/4}$. This implies

$$(26) \quad \left| \sum_{i=1}^{N^{1-\alpha'}} (\theta_i - E\theta_i) \right| \geq |N^{(1-\alpha')/2} - N^{1-\alpha'} \cdot N^{(1-\alpha')/4}| \geq \frac{N^{3(1-\alpha')/4}}{2},$$

from which we conclude, using Lemma 11 of [1],

$$P_N\left( \sum_{i=1}^{N^{1-\alpha'}} 1_{\{\zeta_i \geq (1-\alpha')\sqrt{2\gamma}(\log N)/2\}} \leq N^{(1-\alpha')/2} \right)$$

$$\leq P_N\left( \left| \sum_{i=1}^{N^{1-\alpha}} (\theta_i - E\theta_i) \right| \geq \frac{N^{3(1-\alpha')/4}}{2} \right)$$

$$\leq \exp\left( -\frac{N^{3(1-\alpha')/2}}{4(2N^{1-\alpha'} + N^{3(1-\alpha')/4})/3} \right)$$

$$\leq \exp(-cN^{(1-\alpha')/2}).$$

By (23) and (24), this is more than we need to prove (15).   □

PROOF OF THEOREM 3.1.   Fix $\alpha \in (1/2, 1)$. From the previous lemma we know that we can find some $\kappa = \kappa(\alpha) > 0$, such that we can assume that at least $N^\kappa$ of the $\varphi_B, B \in \Pi_\alpha$, are positive. We use the notation of the previous section, and define, for $1 \leq k \leq K + 1$ and $\varepsilon > 0$, the event

$$
\begin{aligned}
A_k &:= A_k(\varepsilon, \alpha, K) \\
&= \bigcup_{B' \in \Gamma_{\alpha_k}} \bigcup_{B \in \Gamma_{B', \alpha_{k+1}}} \Big\{ |\varphi_{B'} - E_N(\varphi_B | \mathcal{F}_{\alpha_k})| \\
&\qquad\qquad\qquad\qquad \geq \varepsilon \lambda \alpha 2 \sqrt{2\gamma} \frac{1}{K} \Big( 1 - \frac{1}{K} \Big) \log N \Big\}.
\end{aligned}
$$

By Lemma 2.11, $\mathrm{var}(\varphi_{B'} - E(\varphi_B | \mathcal{F}_{\alpha_k}) | \mathcal{F}_{\alpha_{k+1}}) \leq c$, and we can bound

$$
\begin{aligned}
P(A_k) &\leq |\Gamma_{\alpha_k}| |\Gamma_{B', \alpha_{k+1}}| \exp\Big( -\frac{\varepsilon^2 \lambda^2 \alpha^2 8 \gamma (1/K^2)(1 - 1/K)^2 (\log N)^2}{2c} \Big) \\
&\leq \exp(-c(\log N)^2). \tag{27}
\end{aligned}
$$

We will later choose $K \geq \varepsilon \lambda$, such that $c$ is independent of $\varepsilon$ and $\lambda$.

On $\bigcap_k A_k^c$, we can apply the tree-argument of [1]. For $k \leq K$, we denote by $\underline{B}^{(k)}$ a sequence of $k$ boxes $B_1 \supset B_2 \supset \cdots \supset B_k$, where $B_i \in \Gamma_{\alpha_i}, 1 \leq i \leq k$. Set

$$
D_k := \{ \underline{B}^{(k)} : \varphi_{B_i} \geq (\alpha - \alpha_i) 2 \lambda \sqrt{2\gamma} (1 - 1/K) \log N, 1 \leq i \leq k \}.
$$

We show that if on the $k$th scale there are many such sequences, so there will be on the $(k+1)$st scale. Let $n_k := N^{\kappa + 4\alpha(k-1)(1/K)(1-\lambda^2)}$, where $\kappa$ is the same constant as in Lemma 3.2, and define

$$
C_k := \{ |D_k| \geq n_k \}.
$$

Assume that we are on $C_k$. Choose $n_k$ sequences $\underline{B}_j^{(k)} = \{ B_{j,1}, B_{j,2}, \ldots, B_{j,k} \}$, $1 \leq j \leq n_k$ in $D_k$. Note that $B_{j,k} \neq B_{i,k}$ if $i \neq j$, since otherwise the sequences would coincide. Set

$$
\zeta_j := \frac{1}{|\Gamma_{B_{j,k}, \alpha_{k+1}}|} \sum_{B \in \Gamma_{B_{j,k}, \alpha_{k+1}}} 1_{\{\varphi_B - \varphi_{B_{j,k}} \geq \lambda \alpha 2 \sqrt{2\gamma} (1/K)(1 - 1/K) \log N\}}
$$

Note that $|\Gamma_{B_{j,k}, \alpha_{k+1}}| = (N^{\alpha/K}/2)^4$ and, therefore,

$$
C_k \cap C_{k+1}^c \subset C_k \cap \Big\{ \sum_{j=1}^{n_k} \zeta_j \leq n_{k+1} \cdot \frac{16}{N^{4\alpha/K}} \Big\}.
$$

If we set

$$\tilde{\zeta}_j := \frac{1}{|\Gamma_{B_{j,k},\alpha_{k+1}}|} \sum_{B \in \Gamma_{B_{j,k},\alpha_{k+1}}} 1_{\{\varphi_B - E(\varphi_B|\mathcal{F}_{\alpha_k}) \geq (1+\varepsilon)\lambda\alpha 2\sqrt{2\gamma}(1/K)(1-1/K)\log N\}},$$

we have $\zeta_j \geq \tilde{\zeta}_j$ on $A_k^c$ and, therefore,

$$P_N(C_k \cap C_{k+1}^c \cap A_k^c) \leq P_N\left(\sum_{j=1}^{n_k} \tilde{\zeta}_j \leq n_{k+1} \cdot \frac{16}{N^{4\alpha/K}}\right).$$

To bound this probability, we need some large deviation estimates on the binomial variables $\sum_{j=1}^{n_k} \tilde{\zeta}_j$. Note that, due to (14), the $\varphi_B - E_N(\varphi_B|\mathcal{F}_{\alpha_k})$ are centered Gaussian variables with variance

$$\text{var}(\varphi_B|\mathcal{F}_{\alpha_k}) \geq \frac{\alpha}{K}\gamma\log N + c.$$

Therefore,

$$E_N(\tilde{\zeta}_j|\mathcal{F}_{\alpha_k})$$
$$\geq \inf_B P_N\left(\varphi_B - E_N(\varphi_B|\mathcal{F}_{\alpha_k}) \geq (1+\varepsilon)\lambda\alpha 2\sqrt{2\gamma}\frac{1}{K}\left(1-\frac{1}{K}\right)\log N \Big| \mathcal{F}_{\alpha_k}\right)$$
$$\geq \exp\left(-\frac{(1+\varepsilon)^2\lambda^2\alpha^2 8\gamma(1/K^2)(1-1/K)^2(\log N)^2}{2\alpha(1/K)\gamma\log N}\right)$$
$$= N^{-4(\alpha/K)\lambda^2(1-1/K)^2(1+\varepsilon)^2}.$$

Thus, on $C_k \cap A_k^c$,

$$C_{k+1}^c \subset \left\{\sum_{j=1}^{n_k}(\tilde{\zeta}_j - E(\tilde{\zeta}_j|\mathcal{F}_{\alpha_k})) \leq n_{k+1}(16/N^{4(\alpha/K)})\right.$$
$$\left. -n_k N^{-4(\alpha/K)\lambda^2(1-1/K)^2(1+\varepsilon)^2}\right\}$$
$$\subset \left\{\left|\sum_{j=1}^{N^\kappa}(\tilde{\zeta}_j - E(\tilde{\zeta}_j|\mathcal{F}_{\alpha_k}))\right| \geq \frac{1}{2}N^{\kappa-4(\alpha/K)\lambda^2(1-1/K)^2(1+\varepsilon)^2}\right\},$$

if, for the last line, $\varepsilon$ is chosen such that $(1-1/K)(1+\varepsilon) < 1$, making the second term dominate (recall $\lambda < 1$). Then Lemma 11 of [1] yields on $C_k \cap A_k^c$

$$P_N(C_{k+1}^c|\mathcal{F}_{\alpha_k}) \leq 2\exp\left(-\frac{N^{2\kappa-8\lambda^2(\alpha/K)(1-1/K)^2(1+\varepsilon)^2}}{2N^\kappa + (2/3)N^{\kappa-4\lambda^2(\alpha/K)(1-1/K)^2(1+\varepsilon)^2}}\right)$$
(28)
$$\leq \exp(-N^{\kappa-8\lambda^2(\alpha/K)(1-1/K)^2(1+\varepsilon)^2}).$$

Choose $K$ large enough, such that $\kappa - \frac{8\alpha}{K} > 0$. Note that $n_1 = N^\kappa$. This implies, using Lemma 3.2 and (28),

$$
P_N(C_K^c) \leq P_N(C_1^c) + \sum_{k=2}^{K} (P_N(C_k^c \cap C_{k-1} \cap A_{k-1}^c) + P_N(A_{k-1}))
$$

$$
= P_N(C_1^c) + \sum_{k=2}^{K} (E_N(P_N(C_k^c | \mathcal{F}_{\alpha_k}) 1_{C_{k-1} \cap A_{k-1}^c}) + P_N(A_{k-1}))
$$

$$
(29) \qquad \leq \exp(-c_1(\log N)^2) + (K) \exp(-N^{\kappa - 8\lambda^2(\alpha/K)(1-1/K)^2(1+\varepsilon)^2})
$$

$$
+ \exp(-c_2(\log N)^2)
$$

$$
\leq \exp(-c(\log N)^2).
$$

Let now $\mathcal{H}_N(a) := \{x \in V_N^\delta : \varphi_x \geq 2\sqrt{2\gamma} a \log N\}$. We consider the event

$$
L_K = L_K(\alpha, \lambda) := \{|\mathcal{H}_N(\lambda(\alpha - \alpha_{K-1}))| \leq n_{K-1}\}.
$$

Note that

$$
P_N(L_K) \leq P_N(|\{B \in \Pi_{\alpha_K} : \varphi_{x_B} \geq 2\sqrt{2\gamma}\lambda(\alpha - \alpha_{K-1})\log N\}| \leq n_{K-1}).
$$

This implies

$$
P(L_K \cap C_K) \leq E_N(P(|\{B \in \Pi_{\alpha_K} : \varphi_{x_B} \geq 2\sqrt{2\gamma}\lambda(\alpha - \alpha_{K-1})\log N\}|
$$

$$
\leq n_{K-1} | \mathcal{F}_{\alpha_K}) 1_{C_K}).
$$

On $C_K \cap L_K$ we have at least $n_K$ boxes $B \in \Pi_{\alpha_K}$ with $\varphi_B \geq 2\sqrt{2\gamma}\lambda(\alpha - \alpha_K)\log N$, and only for at most $n_{K-1}$ of them we have $\varphi_{x_B} \geq 2\sqrt{2\gamma}\lambda(\alpha - \alpha_K)\log N$. Thus, for at least $n_K - n_{K-1}$ boxes, $\varphi_{x_B} - \varphi_B \leq \mu_K \log N$, with $\mu_K := 2\sqrt{2\gamma}\lambda(\alpha_K - \alpha_{K-1})$. Now we use the fact that, conditional on $\mathcal{F}_{\alpha_K}$, the $\varphi_{x_B} - \varphi_B$ are independent centered Gaussian with variance equal to $\gamma\alpha_K \log N$, and that $\alpha_K - \alpha_{K-1} = -\frac{\alpha}{K} < 0$, and $n_{K-1} = n_K N^{-(4\alpha/K)(1-\lambda^2)}$ to obtain

$$
P_N(|\{B \in \Pi_{\alpha_K} : \varphi_{x_B} \geq 2\sqrt{2\gamma}\lambda(\alpha - \alpha_{K-1})\log N\}| \leq n_{K-1} | \mathcal{F}_{\alpha_K})
$$

$$
\leq P_N(|\{B \in \Pi_{\alpha_K} : \varphi_{x_B} - \varphi_B \leq \mu_K \log N\}| \geq n_K - n_{K-1} | \mathcal{F}_{\alpha_K})
$$

$$
(30) \qquad \leq P_N\left(\varphi_{x_B} - \varphi_B \leq -\frac{\alpha}{K} 2\sqrt{2\gamma}\lambda(1 - 1/K)\log N | \mathcal{F}_{\alpha_K}\right)^{(1 - N^{-(4\alpha/K)(1-\lambda^2)})n_K}
$$

$$
\leq \exp\left(-4\lambda^2 \frac{\alpha}{K}(1 - 1/K)^2(\log N)(1 - N^{-(4\alpha/K)(1-\lambda^2)})n_K\right)
$$

$$
\leq \exp(-c(\log N)^2).
$$

To complete the proof, we get from (29) and (30), using $\alpha - \alpha_{K-1} = \alpha(1 - \frac{2}{K})$,

$$P_N\left(\left|\mathcal{H}_N\left(\lambda\alpha\left(1 - \frac{2}{K}\right)\right)\right| \leq n_{K-1}\right) \leq P_N(L_K \cap C_K) + P_N(C_K^c)$$

(31)
$$\leq \exp(-c(\log N)^2).$$

We can now choose $K$ large enough and $\alpha$ close to 1, such that with (31)

$$P_N(|\{x \in V_N^\delta : \varphi_x \geq 2\sqrt{2\gamma}\lambda \log N\}| \leq N^{4(1-\lambda^2)-\varepsilon})$$
$$\leq P_N\left(\left|\mathcal{H}_N\left(\lambda\alpha\left(1 - \frac{2}{K}\right)\right)\right| \leq n_{K-1}\right)$$
$$\leq \exp(-c(\log N)^2). \qquad \square$$

**4. Probability to stay positive.** Having obtained the same result for the maximum of the interface as in the case of the 2-dimensional lattice free field, we can again use the strategy of [1].

PROOF OF THEOREM 1.3, THE LOWER BOUND. First, note that by a density argument, $\mathcal{C}_V^2(D) = \inf\{\frac{1}{2}\int_V |\Delta h|^2 \, dx : h \in C_0^\infty(V), h \geq 1 \text{ a.e. on } D\}$, where $C_0^\infty(V)$ denotes the infinitely often differentiable functions on $V$ which vanish at $\partial V$. Choose a function $f \in C_0^\infty(V), f \geq 0, f = 1$ on $D$, and a number $a > 2\sqrt{2\gamma}$. Set $\tilde{\varphi}_x := \varphi_x + a \log N f(\frac{x}{N})$. Then $\{\tilde{\varphi}_x\}_{x \in V_N}$ is a Gaussian family with covariances $G_N(x, y), x, y \in V_N$, and expectation $a \log N f(\frac{x}{N})$. Denote the law of this family by $P_N^a$, and let $f_N(x) := f(x/N)$. The relative entropy of $P_N^a$ with respect to $P_N$ is defined as $H_N(P_N^a|P_N) := E_N^a(\log \frac{dP_N^a}{dP_N})$. Note that

$$\frac{dP_N^a}{dP_N}(\varphi) = \exp\left[\frac{1}{2}(\langle \varphi, G_N^{-1}\varphi\rangle_{V_N} - \langle \varphi - a \log N f_N, G_N^{-1}(\varphi - a \log N f_N)\rangle_{V_N})\right],$$

where $\langle \cdot, \cdot \rangle_{V_N}$ denotes the $L_2$-scalar product on $V_N$ and, therefore,

$$E_N^a\left(\log \frac{dP_N^a}{dP_N}\right) = \frac{a^2}{2}(\log N)^2\langle \Delta_N f_N, \Delta_N f_N\rangle_{V_N},$$

from which we conclude

$$\lim_{N \to \infty} \frac{1}{(\log N)^2}H_N(P_N^a|P_N) = \frac{a^2}{2}\|\Delta f\|_{L^2(V)}^2.$$

Moreover,

$$P_N^a((\Omega_N^+)^c) \leq \sum_{x \in D_N} P_N^a(\varphi_x < 0) = \sum_{x \in D_N} P_N(\varphi_x < -a \log N)$$
$$\leq N^4 \exp\left(\frac{-a^2(\log N)^2}{2\gamma \log N}\right) = o(1)$$

as $N \to \infty$. Using the entropy inequality (see, e.g., [7], Appendix B.3), we have

$$\log \frac{P_N(\Omega_N^+)}{P_N^a(\Omega_N^+)} \geq -\frac{H_N(P_N^a|P_N) + e^{-1}}{P_N^a(\Omega_N^+)}$$

and, hence,

$$\liminf_{N \to \infty} \frac{1}{(\log N)^2} \log P_N(\Omega_N^+) \geq -\frac{a^2}{2} \|\Delta f\|_{L^2(V)}$$

for any choice of $a$ and $f$ as above. Optimizing over $a$ and $f$ gives the lower bound. $\square$

PROOF OF THE UPPER BOUND.    Fix $\beta > 0$. For $K \in \mathbb{N}, \alpha \in (1/2, 1)$ define

$$E_{K,\beta,\alpha} := \{\sharp\{B \in \Pi_\alpha : B \subset D_N, \varphi_B \leq (2\sqrt{2\gamma} - \beta) \log N\} \leq K\},$$

the event that we have few boxes $B \in \Pi_\alpha$ with $\varphi_B \leq (2\sqrt{2\gamma} - \beta) \log N$. We will now show that the probability that $\Omega_N^+$ occurs on $E_{K,\beta,\alpha}^c$ is small. If $\eta > 0, \varepsilon \in (0, 1/2), \alpha \in (0, 1)$, let

$$A := \bigcup_{B \in \Pi_\alpha} \bigcup_{x \in B^{(\varepsilon)}} \{|\varphi_B - E_N(\varphi_x|\mathcal{F}_\alpha)| \geq \eta \log N\},$$

where $B^{(\varepsilon)}$ is the set of points $x \in B$, which are contained inside a box of side-length $\varepsilon N^\alpha$ and center $x_B$. We split

$$P_N(E_{K,\beta,\alpha}^c \cap \Omega_{D_N}^+) \leq E_N(P_N((E_{K,\beta,\alpha}^c \cap \Omega_{D_N}^+)|\mathcal{F}_\alpha)1_{A^c}) + P_N(A).$$

But, by Lemma 2.11, we find

$$P_N(A) \leq N^4 \exp\left(-\frac{\eta^2(\log N)^2}{c\varepsilon}\right) \leq \exp\left(-\frac{c'\eta^2(\log N)^2}{\varepsilon}\right).$$

We can choose $\varepsilon$ arbitrarily small; our choice will be such that $\frac{c'\eta^2}{\varepsilon} \geq 8\gamma\mathcal{C}_V^2(D) + 1$. Fix $B \in \Pi_\alpha$, and set $B^{(\varepsilon)} := \{x \in B : \mathrm{dist}(x, \partial B) \geq \varepsilon N^\alpha\}$. The idea is to apply Theorem 1.2 to the field $(\varphi_x - E_N(\varphi_x|\mathcal{F}_\alpha))_{x \in B}$ conditional on $\mathcal{F}_\alpha$. We get

$$P_N\left(\sup_{x \in B^{(\varepsilon)}} (\varphi_x - E_N(\varphi_x|\mathcal{F}_\alpha)) \leq (2\sqrt{2\gamma} - \beta) \log N|\mathcal{F}_\alpha\right)$$

$$\leq P_N\left(\sup_{x \in B^{(\varepsilon)}} (\varphi_x - E(\varphi_x|\mathcal{F}_\alpha)) \leq (2\sqrt{2\gamma} - \beta/2) \log N^\alpha|\mathcal{F}_\alpha\right)$$

$$\leq \exp(-c(\log N)^2),$$

where $c = c(\varepsilon, \beta)$ if $\alpha \in (\alpha_0(\beta), 1)$ for some $\alpha_0(\beta) > 0$. Therefore, on $A^c \cap \{\varphi : \varphi_B \leq (2\sqrt{2\gamma} - \beta) \log N\}$ we have, if $\eta \leq \beta/2$,

$$P_N\left(\inf_{x \in B} \varphi_x \geq 0 \Big| \mathcal{F}_\alpha\right)$$
$$\leq P_N\left(\inf_{x \in B^{(\varepsilon)}} (\varphi_x - E_N(\varphi_x | \mathcal{F}_\alpha)) \geq -(2\sqrt{2\gamma} - \beta/2) \log N | \mathcal{F}_\alpha\right)$$
$$\leq \exp(-c(\log N)^2)$$

if $\alpha \geq a_0(\beta)$. This implies

$$P_N(E_{K,\beta,\alpha}^c \cap \Omega_N^+) \leq \binom{N^{4-4\alpha}}{K} (\exp(-c(\log N)^2))^K$$
$$+ \exp(-(8\gamma \mathcal{C}_V^2(D) + 1)(\log N)^2)$$
$$\leq \exp((4 - 4\alpha)K \log N - cK(\log N)^2)$$
$$+ \exp(-(8\gamma \mathcal{C}_V^2(D) + 1)(\log N)^2)$$
$$\leq \exp(-(8\gamma \mathcal{C}_V^2(D) + 1)(\log N)^2)$$

if we choose $K$ large enough such that $cK/2 \geq 8\gamma \mathcal{C}_V^2(D) + 1$.

This means we now only need to consider $E_{K,\beta,\alpha} \cap \Omega_{D_N}^+$. In this case, for any function $f \geq 0, f \in C^2(D)$, we have

$$\frac{1}{|\Pi_\alpha|} \sum_{B \in \Pi_\alpha, B \subset D_N} f(x_B/N)\varphi_B$$
$$\geq (2\sqrt{2\gamma} - \beta) \log N \left(\frac{1}{|\Pi_\alpha|} \sum_{B \in \Pi_\alpha, B \subset D_N} f(x_B/N) - \frac{K\|f\|_\infty}{|\Pi_\alpha|}\right).$$

Therefore,

$$P_N(E_{K,\beta,\alpha} \cap \Omega_{D_N}^+)$$
$$\leq \exp\left(-\frac{((2\sqrt{2\gamma} - \beta) \log N(1/|\Pi_\alpha| \sum_B f(x_B/N) - cN^{-4(1-\alpha)}))^2}{2 \operatorname{var}_N(1/|\Pi_\alpha| \sum_B f(x_B/N)\varphi_B)}\right).$$

Applying Lemmas C.1 and C.2 completes the proof. $\square$

## 5. Entropic repulsion.

Here we need to use a different approach than in the lattice free field case, since the FKG property does not hold.

PROOF OF PROPOSITION 1.4.   Let $P_N^+(\cdot) := P_N(\cdot | \Omega_N^+)$. We use the notation of Section 3, in particular, the box-structure, and first assume $x = 0$. Set $\overline{\varphi}_{\varepsilon N} := \overline{\varphi}_{\varepsilon N}(x)$. We claim that, on the set $\{\overline{\varphi}_{\varepsilon N} \le (2\sqrt{2\gamma} - \eta) \log N\} \cap \Omega_N^+$, there exists $\delta > 0$ such that

$$\sharp\{x \in V_{\varepsilon N} : \varphi_x \le (2\sqrt{2\gamma} - \eta/2) \log N\} \ge \delta |V_{\varepsilon N}|.$$

If this was not the case, we would have

$$(1-\delta)(2\sqrt{2\gamma} - \eta/2) \log N \le \overline{\varphi}_{\varepsilon N} \le (2\sqrt{2\gamma} - \eta) \log N,$$

which is impossible if $\delta$ is small enough such that $(1-\delta)(2\sqrt{2\gamma} - \eta/2) > (2\sqrt{2\gamma} - \eta)$. Therefore, if $\alpha \in (0,1)$, there exists a shift of the $N^\alpha$-sublattice $\Pi_\alpha$ such that, for this particular shift,

$$P_N^+(\sharp\{x \in V_{\varepsilon N} : \varphi_x \le (2\sqrt{2\gamma} - \eta/2) \log N\} \ge \delta |V_{\varepsilon N}|)$$

$$= P_N^+\left(\frac{1}{|V_{\varepsilon N}|} \sum_{x \in V_{\varepsilon N}} 1_{\{\varphi_x \le (2\sqrt{2\gamma} - \eta/2) \log N\}} \ge \delta\right)$$

$$\le P_N^+\left(\frac{1}{|\{B \in \Pi_\alpha, x_B \in V_{\varepsilon N}\}|} \sum_{B \in \Pi_\alpha, x_B \in V_{\varepsilon N}} 1_{\{\varphi_{x_B} \le (2\sqrt{2\gamma} - \eta/2) \log N\}} \ge \delta\right).$$

(This is true since $\frac{1}{|V_{\varepsilon N}|} \sum_{x \in V_{\varepsilon N}} 1_{\{\varphi_x \le (2\sqrt{2\gamma} - \eta/2) \log N\}}$ is the average over all possible such shifts of the $N^\alpha$-lattice.) Let $S_\alpha := \{B \in \Pi_\alpha, x_B \in V_{\varepsilon N}\}$ for this particular $\Pi_\alpha$. Choose $0 < \delta' < \delta$. Then

$$P_N^+\left(\frac{1}{|S_\alpha|} \sum_{B \in S_\alpha} 1_{\{\varphi_{x_B} \le (2\sqrt{2\gamma} - \eta/2) \log N\}} \ge \delta\right)$$

(32)
$$\le P_N^+\left(\frac{1}{|S_\alpha|} \sum_{B \in S_\alpha} 1_{\{\varphi_B \le (2\sqrt{2\gamma} - \eta/4) \log N\}} \ge \delta'\right)$$

$$+ P_N^+\left(\frac{1}{|S_\alpha|} \sum_{B \in S_\alpha} 1_{\{\varphi_B - \varphi_{x_B} > (\eta/4) \log N\}} \ge (\delta - \delta')\right).$$

We have $|S_\alpha| \ge c\varepsilon N^{4(1-\alpha)}$. Thus,

$$P_N^+\left(\frac{1}{|S_\alpha|} \sum_{B \in S_\alpha} 1_{\{\varphi_B \le (2\sqrt{2\gamma} - \eta/4) \log N\}} \ge \delta'\right)$$

$$\le P_N^+(\sharp\{B \in \Pi_\alpha : \varphi_B \le (2\sqrt{2\gamma} - \eta/4) \log N\} \ge c\delta'\varepsilon N^{4(1-\alpha)}).$$

But in the proof of the upper bound of Theorem 1.3 we have seen that

$$P_N(E_{k,\beta,\alpha}^c \cap \Omega_N^+) \le \exp(-(8\gamma \mathcal{C}_V^2(D) + 1)(\log N)^2),$$

hence, for large enough $N$,

$$P_N^+(\sharp\{B \in \Pi_\alpha : \varphi_B \le (2\sqrt{2\gamma} - \eta/4) \log N\} \ge c\delta'\varepsilon N^{4(1-\alpha)})$$

$$\le \exp(-c(\log N)^2).$$

Thus, what is left is the second term in (32). Note

$$P_N(\varphi_B - \varphi_{x_B} > (\eta/4)\log N|\mathcal{F}_\alpha) \le \exp(-c\eta^2\log N).$$

Let $\theta_B := 1_{\{\varphi_B - \varphi_{x_B} > (\eta/4)\log N\}}$. As in the proof of Theorem 1.2, we have, using Lemma 11 of [1], for large $N$,

$$P_N\left(\sum_{B \in S_\alpha} 1_{\{\varphi_B - \varphi_{x_B} > (\eta/4)\log N\}} \ge (\delta - \delta')|S_\alpha|\right)$$

$$\le P_N\left(\left|\sum_{B \in S_\alpha}(\theta_B - E\theta_B)\right| \ge c\varepsilon N^{4(1-\alpha)}((\delta - \delta') - N^{-c'\eta^2})\right)$$

$$\le P_N\left(\left|\sum_{B \in S_\alpha}(\theta_B - E\theta_B)\right| \ge c\varepsilon(\delta - \delta')N^{4(1-\alpha)}\right)$$

$$\le 2\exp(-c\varepsilon(\delta - \delta')N^{4(1-\alpha)}).$$

Together with Theorem 1.3, this proves

$$\lim_{N \to \infty} P_N(\overline{\varphi}_{\varepsilon N} \le (2\sqrt{2\gamma} - \eta)\log N \,|\, \Omega_N^+) = 0$$

if $x = 0$. For arbitrary $x$ repeat the argument with a shifted grid.   □

## APPENDIX A: NORM ESTIMATES

In this section we prove some basic estimates on the discrete Sobolev norms which are used in the proof of the regularity for the solution of the Dirichlet problem. Recall

$$E_1 = \{v : V_N \cup \partial_2 V_N \to \mathbb{R} : v(x) = 0 \ \forall x \in \partial_2 V_N\}$$

and for $v, w \in E_1$ from Section 2,

$$\mathcal{D}(v, w) := \sum_{x \in V_N} \Delta v(x)\Delta w(x) + \sum_{x \in \partial_- V_N} r(x)v(x)w(x).$$

Note that the notation $\mathcal{D}(v, w)$ and $E_1$ depend on $N$. We identify $v \in E_1$ with the function we obtain if we extend $v$ to all of $\mathbb{Z}^d$ by setting it equal to 0 on the whole of $V_N^c$.

Lemma A.1.   *Let $v \in E_1$. There exists a constant $c$ depending on the dimension such that*

$$\sum_{x \in V_{N+1}} \sum_{i=1}^{d} \sum_{j=1}^{d} (\nabla_i \nabla_j v(x))^2 \le c\mathcal{D}(v, v).$$

PROOF.   Expanding the square gives

$$(2d)^2 \sum_{x \in V_N} (\Delta v(x))^2$$

$$= \sum_{x \in V_N} \sum_{i,j=1}^d (4v(x)^2 - 2v(x)v(x+e_i) - 2v(x)v(x-e_i)$$

$$(33) \qquad\qquad - 2v(x)v(x+e_j) - 2v(x)v(x-e_j)$$

$$+ v(x+e_i)v(x+e_j) + v(x+e_i)v(x-e_j)$$

$$+ v(x-e_i)v(x+e_j) + v(x-e_i)v(x-e_j)).$$

Now, taking the geometry of $V_N$ and the 0-boundary conditions outside $V_N$ into consideration, we can shift the summation, and obtain for any $e_i$ with $|e_i| = 1$,

$$\sum_{x \in V_N} v(x)^2 = \sum_{x \in V_{N+1}} v(x)^2 = \sum_{x \in V_{N+1}} v(x+e_i)^2$$

$$= \sum_{x \in V_{N+1}} v(x+e_i+e_j)^2 + \sum_{\substack{x \notin V_{N+1}: \\ x+e_i+e_j \in V_N}} v(x+e_i+e_j)^2.$$

Similarly, we have

$$\sum_{x \in V_N} v(x)v(x-e_i)$$

$$= \sum_{x \in V_N} v(x+e_i)v(x)$$

$$= \sum_{x \in V_{N+1}} v(x+e_i+e_j)v(x+e_j) + \sum_{\substack{x \notin V_{N+1}: \\ x+e_i+e_j \in V_N \\ x+e_j \in V_N}} v(x+e_i+e_j)v(x+e_j)$$

and

$$\sum_{x \in V_N} v(x-e_i)v(x+e_j) = \sum_{x \in V_N} v(x+e_i+e_j)v(x) = \sum_{x \in V_{N+1}} v(x+e_i+e_j)v(x).$$

Furthermore, if $i \neq j$,

$$\sum_{x \in V_N} v(x-e_i)v(x-e_j) = \sum_{x \in V_N} v(x+e_i)v(x+e_j) = \sum_{x \in V_{N+1}} v(x+e_i)v(x+e_j)$$

and

$$\sum_{x \in V_N} v(x-e_i)^2 = \sum_{x \in V_{N+1}} v(x+e_i)^2 - \sum_{\substack{x \in V_N: \\ x+e_i \notin V_N}} v(x)^2$$

and, finally,

$$\sum_{x \in V_N} v(x+e_i)^2 = \sum_{x \in V_{N+1}} v(x+e_i)^2 - \sum_{\substack{x \in V_N: \\ x-e_i \notin V_N}} v(x)^2.$$

We define the following quantities:

$$T_1 := \sum_{i,j=1}^{d} \sum_{x \notin V_{N+1}} v(x+e_i+e_j)^2 \geq 0,$$

$$T_2 := \sum_{i=1}^{d} \sum_{\substack{x \in V_N: \\ x+e_i \notin V_N}} v(x)^2 \quad \text{and} \quad T_3 := \sum_{i=1}^{d} \sum_{\substack{x \in V_N: \\ x-e_i \notin V_N}} v(x)^2.$$

Note $T_2 + T_3 \leq \sum_{x \in \partial_- V_N} r(x)v(x)^2$. By the above considerations, the right-hand side of (33) can be rewritten and bounded as follows:

$$(2d)^2 \sum_{x \in V_N} (\Delta v(x))^2$$

$$= \sum_{x \in V_N} \sum_{i,j=1}^{d} (v(x)^2 + v(x+e_i)^2 + v(x+e_j)^2 + v(x+e_i+e_j)^2$$

$$- 2v(x)v(x+e_i) - 2v(x+e_i+e_j)v(x+e_j) - 2v(x)v(x+e_j)$$

$$- 2v(x+e_j)v(x) + v(x+e_i)v(x+e_j) + 2v(x+e_i+e_j)v(x)$$

$$+ v(x+e_i+e_j)v(x+e_i))$$

$$+ T_1 - T_2 - T_3$$

$$\geq \sum_{i,j=1}^{d} \sum_{x \in V_{N+1}} (\nabla_i \nabla_j v(x))^2 - \sum_{x \in \partial_- V_N} r(x)v(x)^2.$$

Thus,

$$\sum_{i,j=1}^{d} \sum_{x \in V_{N+1}} (\nabla_i \nabla_j v(x))^2 \leq (2d)^2 \sum_{x \in V_N} (\Delta v(x))^2 + \sum_{x \in \partial_- V_N} r(x)v(x)^2$$

$$\leq (2d)^2 \mathcal{D}(v,v),$$

which proves the lemma.  $\square$

LEMMA A.2.  *Let $v \in E_1$. There exists $c > 0$ such that*

$$\sum_{x \in V_N} v(x)^2 \leq cN^2 \left( \sum_{x \in V_N} \sum_{i=1}^{d} (\nabla_i v(x))^2 + \sum_{x \in \partial_- V_N} r(x)v(x)^2 \right).$$

PROOF.  Let $x \in V_N$ and denote $A_x^i := \{y \in V_N : \exists k \in \mathbb{Z}$ such that $y = x + k \cdot e_i\}$. Then

$$v(x)^2 = (v(x) - v(x + e_i) + v(x + e_i) - v(x + 2e_i) + \cdots + v(x + k_0 e_i))^2,$$

where $k_0 \in \mathbb{N}$ such that $x + k_0 e_i \in \partial_- V_N$. Obviously $k_0 \leq 2N$, thus, using the fact that $(a + b)^2 \leq 2a^2 + 2b^2$ for real numbers $a, b$, we get

$$v(x)^2 \leq 2N((v(x) - v(x + e_i))^2 + \cdots$$
$$+ (v(x + (k_0 - 1)e_i) - v(x + k_0 e_i))^2 + v(x + k_0 e_i)^2).$$

In the same way, we obtain

$$v(x)^2 \leq 2N((v(x) - v(x - e_i))^2 + \cdots + v(x + k_1 e_i)^2)$$

for some $k_1 \leq 2N$, with $x - k_1 e_1 \in \partial_- V_N$. This gives

$$\sum_{x \in V_N} v(x)^2 \leq 2 \sum_{x \in V_N} N \left( \sum_{y \in A_x^i} (v(y) - v(y + e_i))^2 + \sum_{y \in \partial_- V_N \cap A_x^i} v(y)^2 \right)$$
$$\leq cN^2 \left( \sum_{x \in V_N} (v(x) - v(x + e_i))^2 + \sum_{x \in \partial_- V_N} r(x)v(x)^2 \right).$$

Since this inequality holds for any $1 \leq i \leq d$, the lemma is proven.  □

LEMMA A.3.  *Let $v \in E_1$. There exists $c > 0$ such that, for all $1 \leq i \leq d$,*

$$\sum_{x \in V_N} (v(x + e_i) - v(x))^2 \leq cN^2 \left( \sum_{x \in V_N} (\nabla_i \nabla_i v(x))^2 + \sum_{x \in \partial_- V_N} r(x)v(x)^2 \right).$$

PROOF.  Let $h(x) := \nabla_i v(x)$ and repeat the arguments of the proof of Lemma A.2.  □

From Lemmas A.2 and A.3 the following is clear:

COROLLARY A.4.  *Let $v \in E_1$. There exists $c > 0$ such that*

$$\|v\|_{H^2(V_N)}^2 \leq cN^4 \left( \sum_{x \in V_N} \sum_{i,j=1}^{d} (\nabla_i \nabla_j v(x))^2 + \sum_{x \in \partial_- V_N} r(x)v(x)^2 \right).$$

REMARK A.5.  Iterating this procedure, one evidently obtains for any $v : V_N \cup \partial_k V_N \to \mathbb{R}$ such that, $v(x) = 0$ for $x \in \partial_k V_N$, that

$$\|v\|_{H^k(V_N)}^2 \leq cN^{2k} \left( \sum_{x \in V_N} \sum_{\alpha : |\alpha| = k} (\nabla^\alpha v(x))^2 + \sum_{x \in \partial_- V_N} r(x)v(x)^2 \right).$$

COROLLARY A.6.    *Let $v \in E_1$.  There is $c > 0$ such that*

$$\|v\|^2_{H^2(V_N)} \leq cN^4 \mathcal{D}(v,v).$$

PROOF.    From Lemma A.1 and Corollary A.4 we obtain

$$\|v\|^2_{H^2(V_N)} \leq \|v\|^2_{H^2(V_{N+1})} \leq c'(N+1)^4 \mathcal{D}(v,v) \leq cN^4 \mathcal{D}(v,v). \qquad \square$$

REMARK A.7.    This also proves that $\mathcal{D}(\cdot,\cdot)$ is positive definite.

## APPENDIX B: DISCRETE SOBOLEV IMBEDDING

The following results are the discrete analogues of the Sobolev Imbedding Theorems. For completeness, we include the proofs of the versions we use.

PROPOSITION B.1.    *Let $f: \mathbb{Z}^d \to \mathbb{R}$ such that $f(x) = 0$ on $V_N^c$, and $\|f\|_{H^k(V_N)} \leq cN^{d/2}$ for some constant $c$ independent of $N$. If $k > d/2$, then there exists $C > 0$ independent of $N$ such that $\sup_{x \in V_N} |f(x)| < C$.*

PROOF.    Let $\widehat{f}(t) = \sum_{x \in \mathbb{Z}^d} f(x) e^{i\langle t, x \rangle}$ denote the Fourier transform of a function $f: \mathbb{Z}^d \to \mathbb{R}$. Then we have

$$\begin{aligned}
\widehat{\nabla_k f}(t) &= \sum_{x \in \mathbb{Z}^d} (f(x+e_k) - f(x)) e^{i\langle t, x \rangle} \\
&= \sum_{x \in \mathbb{Z}^d} (f(x) e^{i\langle t, x-e_k \rangle} - f(x) e^{i\langle t, x \rangle}) \\
&= \widehat{f}(t)(e^{-it_k} - 1).
\end{aligned}$$

Iterating, we obtain

$$(34) \qquad \widehat{\nabla_{k_1} \cdots \nabla_{k_l} f}(t) = \widehat{f}(t)(e^{-it_{k_1}} - 1) \cdots (e^{-it_{k_l}} - 1).$$

By (34), using the Taylor expansion, we have, for any $j \in \mathbb{N}$,

$$|\widehat{f}(t)|^2 \cdot |t|^{2j} \leq c \cdots |\widehat{f}(t)|^2 |(e^{-it_{k_1}} - 1) \, \cdots (e^{-it_{k_l}} - 1)|^2 \leq |\widehat{\nabla_{k_1} \cdots \nabla_{k_l} f}(t)|^2.$$

This yields

$$\begin{aligned}
\int_{[-\pi,\pi]^d} &|\widehat{f}(t)| \, dt \\
&= \int_{V_N} \frac{1}{(1 + N^2|t|^2)^{l/2}} (1 + N^2|t|2)^{l/2} |\widehat{f}(t)| \, dt \\
&\leq \left( \int_{[-\pi,\pi]^d} \frac{1}{(1 + N^2|t|^2)^l} \, dt \right)^{1/2}
\end{aligned}$$

$$\times \left( \int_{[-\pi,\pi]^d} (1 + N^2|t|^2)^l |\widehat{f}(t)|^2 \, dt \right)^{1/2}$$

$$\leq c N^{-l} \cdot \left( \int_{[-\pi,\pi]^d} \sum_{j=0}^{l} (N|t|)^{2j} |\widehat{f}(t)|^2 \, dt \right)^{1/2}$$

$$\leq c N^{-l} \|f\|_{H^l(V_N)} \leq c N^{d/2-l},$$

using the Plancherel Theorem. Thus, we get, by the inverse Fourier transform,

$$|f(x)| = \left| c \int_{[-\pi,\pi]^d} \widehat{f}(t) e^{-i\langle t,x \rangle} \, dt \right| \leq \int_{[-\pi,\pi]^d} |\widehat{f}(t)| \, dt \leq c N^{d/2-l}. \qquad \square$$

This implies the following:

COROLLARY B.2. *Let* $f \colon \mathbb{Z}^d \to \mathbb{R}$ *such that* $f(x) = 0$ *on* $V_N^c$, *and* $\|f\|_{H^k(V_N)} \leq c N^{d/2}$ *for some constant* $c$ *independent of* $N$. *If* $k > d/2 + l$, *then there exists* $C > 0$ *independent of* $N$ *such that* $\sup_{x \in V_N} |\nabla^\alpha f(x)| \leq \frac{C}{N^{|\alpha|}}$ *for all* $0 \leq |\alpha| \leq l$.

## APPENDIX C: COMPUTATION OF THE CONSTANT $\mathcal{C}_V^2(D)$

We still need to show the convergence toward the second-order capacity $\mathcal{C}_V^2(D)$ in the upper bound of Theorem 1.3. This is analogous to a similar statement in the higher-dimensional case; compare [10]. Let

$$H_0^2(V_N) := \{ f \in H^2(V_N) \colon f(x) = 0 \ \forall x \in \partial_- V_N \}$$

and

$$C_0^\infty(V_N) := \{ f \colon V_N \to \mathbb{R} \colon |\nabla^\alpha f| \leq c/N^{|\alpha|}, \alpha \in \mathbb{N}_0^d, f(x) = 0, \forall x \in \partial_- V_N \}.$$

If $f \colon V \to \mathbb{R}$, we write $f_N$ for the function $V_N \to \mathbb{R}$, $f_N(x) := f(x/N)$.

LEMMA C.1.

$$\inf\{ \|\Delta_N h\|_{L_2(V_N)}^2 \colon h \in H_0^2, h \geq 1 \text{ on } D_N \}$$

$$= \sup\left\{ \langle 1_{D_N}, f_N \rangle_{D_N} - \frac{1}{2} \langle f_N G_N f_N \rangle \colon f \in L_2(V_N) \colon f = 0 \text{ on } V_N \setminus D_N \right\}$$

$$= \sup\left\{ \frac{\langle 1_{D_N}, f_N \rangle_{D_N}^2}{2 \langle f_N, G_N f_N \rangle_{D_N}} \colon f \in L_2(V_N) \colon f = 0 \text{ on } V_N \setminus D_N \right\}.$$

Proof. We start with the first equality. Since $E_0(V_N)$ is finite dimensional, there exists a minimizer $h_N^{(0)}$. Obviously, $h_N^{(0)} = 1$ on $D_N$. Furthermore, $\Delta^2 h_N^{(0)} = 0$ outside $D_N$. To see this, set $\psi(\varepsilon) = \sum_{x \in V_N} |\Delta h_N^{(0)}(x) + \varepsilon \varphi(x)|$ for any test function $\varphi \colon V_N \cup \partial_2 V_N \to \mathbb{R}$, with $\varphi(x) = 0$ for all $x \in V_N \setminus D_N$. Then $\frac{d\psi}{d\varepsilon}|_{\varepsilon=0} = 0$, because $h_N^{(0)}$ is a minimizer of the norm. But this implies $\langle \Delta^2 h_N^{(0)}, \varphi \rangle_{V_N} = \langle \Delta h_N^{(0)}, \Delta\varphi \rangle_{V_N} = 0$ for all $\varphi$ as above, and thus the claim. Set

$$f_N = \Delta_N^2 h_N^{(0)}.$$

By the fact that $f_N^{(n)} = 0$ outside $D_N$, summation by parts gives

$$2\langle f_N, h_N^{(0)} \rangle_{D_N} - \langle f_N, G_N f_N \rangle_{D_N} = \sum_{x \in V_N} |\Delta h_N^{(0)}|^2.$$

The above yields

$$\sup\left\{ \langle 1_{D_N}, f_N \rangle_{D_N} - \frac{1}{2} \langle f_N G_N f_N \rangle : f \in L_2(V_N) : f = 0 \text{ on } V_N \setminus D_N \right\}$$

$$\geq 2\langle f_N, h_N^{(0)} \rangle_{D_N} - \langle f_N, G_N f_N \rangle_{D_N}$$

$$= \sum_{x \in V_N} |\Delta_N h_N^{(0)}|^2,$$

which is one direction in the first equation. The other direction is an elementary calculation.

The second equation follows by expanding $f$ in a basis of eigenvectors of the symmetric matrix $G_N$. Maximizing shows that both sides are equal to $\sum_{i \in \mathbb{N}} \frac{\langle e_i, 1_D \rangle^2}{\lambda_i}$, where the $e_i$ are the eigenvectors and $\lambda_i$ the corresponding eigenvalues. □

Lemma C.2. *With the above notation,*

$$\lim_{N \to \infty} \inf\{\|\Delta_N h\|_{L_2(V_N)}^2 : h \in H_0^2, h \geq 1 \text{ on } D_N\} = \mathcal{C}_V^2(D).$$

Proof. $\{h \in H_0^2(V) : h \geq 1_D\}$ is a closed convex subset of the Hilbert space $H_0^2(V)$ and, therefore, there exists a minimizer $h_0$ for $\int_V |\Delta h|^2 \, dx$. For every $n \in \mathbb{N}$, the discretization $h_{0,N}(x) := h_0(x/N)$ belongs to $H_0^2(V_N)$, which proves one direction. Let $\varepsilon > 0$. For every $N \in \mathbb{N}$, we can find $\tilde{h}^{(N)} \in H_0^2(V)$ such that $\tilde{h}^{(N)} \geq 1_D$ and the discretization $\tilde{h}_N^{(N)}$ of $\tilde{h}^{(N)}$ is equal to $h_N^{(0)}$ of the proof of Lemma C.1. If $N$ is large enough, $\|\tilde{h}_N^{(N)}\|_{L_2(V_N)} \geq \|\tilde{h}^{(N)}\|_{L_2(V)} - \varepsilon$. Since $h_0$ is a minimizer, we have $\liminf_{N \to \infty} \|h_N^{(0)}\|_N \geq \liminf_{N \to \infty} \|\tilde{h}^{(N)}\|_{L_2(V)} - \varepsilon \geq \|h_0\|_{L_2(V)} - \varepsilon$. Since $\varepsilon > 0$ was arbitrary, the claim is proven. □

**Acknowledgment.** Many thanks to Erwin Bolthausen for his advice and important discussions.

INSTITUT FÜR MATHEMATIK
WINTERTHURERSTR. 190
CH-8057 ZÜRICH
SWITZERLAND
E-MAIL: noemi.kurt@math.uzh.ch